\documentclass[upref]{amsart}
\usepackage{latexsym,amsmath,amstext,mathrsfs,amssymb,amsthm,amscd,amsopn,verbatim,graphics}

\usepackage{times}
\numberwithin{equation}{section}

\DeclareMathOperator{\id}{id}

\DeclareMathOperator{\conj}{c}

\DeclareMathOperator{\pr}{pr}

\DeclareMathOperator{\ev}{ev}

\DeclareMathOperator{\op}{\mathnormal{op}}

\newtheorem{Thm}{Theorem}[section]
\newtheorem{Prop}[Thm]{Proposition}
\newtheorem{Lem}[Thm]{Lemma}
\newtheorem{Cor}[Thm]{Corollary}

\theoremstyle{remark}
\newtheorem{Rem}[Thm]{Remark}
\newtheorem*{Ack}{Acknowledgment}

\theoremstyle{definition}

\newtheorem{Def}[Thm]{Definition}
\newtheorem{Exa}[Thm]{Example}

\newtheorem*{Reg*}{Regularization procedure}

\newcommand{\calH}{\mathcal{H}}

\newcommand{\calG}{\mathcal{G}}
\newcommand{\calU}{\mathcal{U}}

\newcommand{\tildep}{\widetilde{p}}

\newcommand\qq{\rm}
\newcommand\cmp[1]{{\qq Commun.\ Math.\ Phys.\ \bf #1}}

\newcommand\Kth[1]{{\qq $K$-Theory \bf #1}}

\begin{document} 

\title[The division map of principal\ldots]{The division map of principal bundles with groupoid structure and generalized gauge transformations}
%\date{}

\author[C.~A.~Rossi]{Carlo~A.~Rossi}
\address{Dept.\ of Mathematics---Technion---32000, Haifa---Israel}  
\email{crossi@techunix.technion.ac.il}
\thanks{C.~A.~Rossi acknowledges partial support from the Aly Kaufman Fellowship}

\begin{abstract}
Motivated by the computations done in~\cite{C1}, where I discussed what I called the groupoid of generalized gauge transformations, viewed as a groupoid over the objects of the category $\mathsf{Bun}_{G,M}$ of principal $G$-bundles over a given manifold $M$, I develop in this paper the same arguments for the more general case of {\em principal $\calG$-bundles or principal bundles with structure groupoid $\calG$}, where now $\calG$ is a Lie groupoid.
Most of the concepts introduced in~\cite{C1} can be translated almost verbatim in the framework of principal bundles with structure groupoid $\calG$; in particular, the key r{\^o}le for the construction of generalized gauge transformations is again played by (the equivalent in the framework of principal bundles with groupoid structure of) the division map $\phi_P$.
Moreover, since Hilsum--Skandalis morphisms~\cite{Con},~\cite{HS},\cite{Moer2} are particular principal bundles with structure groupoid, it is possible to develop a notion of Hilsum--Skandalis generalized gauge transformations, by modifying slightly the previous arguments.
\end{abstract}

\maketitle

\tableofcontents

\section{Introduction}
In the paper~\cite{C1}, in order to construct explicit isomorphisms between principal bundles on the space of loops in a manifold $M$, obtained by pulling back a fixed principal bundle $P$ over $M$ w.r.t.\ different kinds of evaluation maps, I introduced the concept of {\em generalized gauge transformation}: namely, as $G$-equivariant (auto)morphisms of a principal $G$-bundle $P$ over $M$ correspond in a bijective way to $G$-equivariant maps from $P$ to $G$, fibre-preserving, $G$-equivariant bundle morphisms between two (a priori) distinct principal $G$-bundles over the same base space $M$ correspond to $G\times G$-equivariant maps from the fibred product bundle of the considered bundles to $G$, which may be viewed as a representation of $G\times G$.
Later, I applied this correspondence to the parallel transport w.r.t.\ a given connection, which can be seen, directly from its well-known $G$-equivariance properties w.r.t.\ initial and final points, as a $G\times G$-equivariant map from the fibred product of $\ev_0^*P$ with $\ev^*P$, where $P$ is a fixed principal $G$-bundle over $M$, and $\ev_0$, resp.\ $\ev$, denotes the evaluation map at the initial point, resp.\ the usual evaluation map; therefore, there exist an explicit bundle isomorphism between $\ev_0^*P$ and $\ev^*P$.
The main tool for establishing the correspondence
\begin{equation}\label{eq-bunmorgengauge}
\begin{aligned}
&\left\{\text{bundle morphisms between $G$-bundles}\right\}\Leftrightarrow\\
&\Leftrightarrow\left\{\text{$G\times G$-equivariant maps from fibred products to $G$} \right\}
\end{aligned}
\end{equation}
is the existence of a canonical map $\phi_P$, attached to any principal $G$-bundle $P$; this map, which is called by MacKenzie~\cite{McK} the {\em division map} of $P$, contains all the informations one needs to characterize the fact that the group $G$ acts on $P$ freely and transitively on each fiber.
Actually, the data of a $G$-invariant surjective submersion from $P$ to $G$, together with the division map $\phi_P$, characterize completely a principal $G$-bundle; the construction of trivializations of $P$ is done explicitly by means of the division map $\phi_P$, see~\cite{McK}.

Let me now write a short story of principal bundles with structure groupoid and of their division map.
The notion of division map for ordinary principal bundles has an analogon in the framework of principal bundles with structure groupoid, i.e.\ when one takes a smooth manifold $P$ (possibly non-Hausdorff), on which a Lie groupoid $\calG$ operates, together with a surjective submersion $\pi$ onto a smooth manifold $M$ (this, in turn, Hausdorff), such that $\pi$ is $\calG$-invariant, and the action of $\calG$ is free and transitive on each fiber of $\pi$.
In fact, the concept of principal bundle with structure groupoid arises naturally in the context of foliations: possibly, the first appearance of this concept was in~\cite{Con}, and later in~\cite{Haef}, where the authors examine the monodromy and holonomy groupoid of a foliated manifold.
Let me notice at this point that both authors prefer to stress the local aspect of principal bundles with groupoid structure, namely they consider mainly a version of nonabelian Cech cohomology for groupoids, and view principal bundles with groupoid structure (or, to be more precise, isomorphisms classes thereof) as Cech cohomology classes on the base space with values in the structure groupoid.
Although this point of view is correct and also, for certain aspects, more fruitful than the one I am going to pursue here, they do not explicitly mention an important piece of the picture, namely the existence of what I call {\em local momenta}; I prefer to skip in this paper any local discussion of principal bundles with structure groupoid, deserving to it a subsequent paper~\cite{C2}. 

In particular, the importance of principal bundles with structure groupoid lies in the notion of {\em generalized morphisms between Lie groupoids} and the strongly related notion of {\em Morita equivalences}: these correspond, roughly, to right principal bundles w.r.t.\ the action of one groupoid, on which another groupoid (a priori distinct) operates from the left in a compatible way, respectively freely, transitively and in a compatible way to the former right action.   
Connes~\cite{Con} also introduced these concepts for Lie groupoids, again from a local point of view, using arguments of nonabelian Cech cohomology; later, Hilsum and Skandalis~\cite{HS} devoted a huge amount of work to generalized morphisms and Morita equivalences.

Approaching later the above subject from a global point of view, M{\oe}rdijk~\cite{Moer1} introduced the notion of division map for a principal bundle with structure groupoid; here, the notion of (global) momentum appeared explicitly in connection to the so-called division map, and the pair formed by (global) momentum and division map is called by M{\oe}rdijk a {\em cocycle on $M$ with values in a Lie groupoid}, where $M$ is the base space of a given principal bundle with structure groupoid.
The notion of cocycle on $M$ with values in a Lie groupoid is equivalent to the global definition of principal bundle with structure groupoid that Mrcun~\cite{Mrcun} adopts for examining in detail the properties of generalized groupoids; this is also illustrated in detail in the book~\cite{Moer2}.

The paper is then organized as follows: Section~\ref{sec-backgr} is simply a review of the main notions concerning Groupoids and Lie Groupoids, the only new thing being (as far as I know) the notion of {\em generalized conjugation in a groupoid}, which is one of the basic notions needed in the rest of the paper.
In Section~\ref{sec-princbun}, I define principal bundles with structure groupoid following~\cite{Mrcun} and~\cite{Moer2}: I will review some basic examples and I will introduce the concept of fibred product of two principal bundles with structure groupoid.
Later, I will introduce the division map of a principal bundle (the terminology is borrowed from the context of ordinary principal bundles, following~\cite{McK}) and I will discuss in detail its properties.
In Section~\ref{sec-eqmorphgengauge}, I will introduce the notion of bundle morphism between principal bundles with structure groupoid and of {\em generalized gauge transformation}; later, using the division map, I will establish the explicit correspondence (\ref{eq-bunmorgengauge}), leading to the notion of {\em groupoid of generalized gauge transformations}.
I devote a small subsection to the invariance property of the division map w.r.t.\ bundle morphisms; this will play a pivotal r{\^o}le in~\cite{C2}, where I plan to discuss in detail the local nature of principal bundles, and hence of generalized morphisms and Morita equivalences between Lie groupoids in the sense specified above.
Finally, Section~\ref{sec-equivhilsumskand} is devoted to the study of {\em Hilsum--Skandalis generalized gauge transformations}: namely, given two groupoids $\calG$ and $\calH$, it is possible to introduce a natural notion of morphism between HS morphisms from $\calG$ to $\calH$, since these are in particular (right) principal $\calH$-bundles.
Therefore, it is sufficient to add one condition to the above notion of bundle morphism to get to the notion of morphisms between HS morphisms; then, by inspecting the properties of the division map of HS morphisms, I derive a notion of HS generalized gauge transformations and I show that morphisms of HS morphisms are in one-to-one correspondence with HS generalized gauge transformations.

I plan to look in the future for possible applications of the correspondence (\ref{eq-bunmorgengauge}) in the framework of gauge theory for principal bundles with structure groupoids: namely, it would be an interesting task to introduce Topological Quantum Field Theories, like e.g.\ Chern--Simons theory, higher-dimensional $BF$-theories in the framework of principal bundles with structure groupoid: in the ordinary case, such constructions rely mainly on notions like principal bundles, associated bundles, connections, etc\ldots
Once one would have introduced and discussed extensively such notions, Correspondence (\ref{eq-bunmorgengauge}) would be a pivotal element in the construction of iterated integrals {\`a} la Chen, representing holonomy, parallel transport and, more generally, borrowing terms from~\cite{CR}, generalized holonomy, which are among the main constructions in Topological Field Theories like Chern--Simons Theory and higher-dimensional $BF$ Theories.
In fact, a concept of connection for principal bundles with structure groupoid was already introduced and briefly discussed in~\cite{Moer2}, where the authors give a ``geometric'' characterization of connection, namely as a particular distribution of the total space of tangent bundle of the principal bundle in question.
It is also possible to give an ``analytic'' characterization of such connections, introducing the concept of connection $1$-forms; this I will do in a forthcoming paper~\cite{C3}.
However, in a forthcoming paper~\cite{C3} I will discuss connections on general principal bundles with structure groupoid by viewing them as generalized gauge transformations between two particular pull-back bundles on the space of curves (parametrized over the unit interval): namely, based on the arguments introduced in the present paper, it is not difficult to see that the expected properties of the parallel transport in the more general framework of principal bundles with structure groupoid fit in into the definition of generalized gauge transformation, with an additional property (which can be regarded as a ``morphism'' property in the context of (quasi) groupoids).
Let me just point out that the theory of connections on a {\em principal $G$-bundle over a Lie groupoid $\Gamma$} has been extensively pursued in a recent paper~\cite{L-GTX}, from where I borrowed the previous notations: namely, $G$ is a Lie group and $\Gamma$ is a Lie groupoid over $\Gamma_0$, the manifold of objects.
Such bundles are, in our context, Hilsum--Skandalis morphisms from $\Gamma$ to $G$, where $G$ can be viewed as a (trivial) Lie groupoid over a point.
However, I would like to pursue in the future the general theory of connections on principal bundles with structure groupoids, and (possibly) to formulate a general Chern--Weil theory, from where it should be possible, in principle, to recover the results of~\cite{L-GTX}.

\begin{Ack}
I thank A.\ S.\ Cattaneo and G.\ Felder for many inspiring suggestions and corrections; I also acknowledge the pleasant atmosphere at the Department of Mathematics of the Technion, where this work was (finally) accomplished.
\end{Ack}

\section{Background definitions: groupoids and Lie groupoids}\label{sec-backgr}
In this section I fix the main notations and conventions regarding the main objects that I consider throughout the paper, namely {\em Lie groupoids}; to begin with, it is better to introduce the concept of a general groupoid.

\begin{Def}\label{def-groupoid}
A {\em groupoid} is a $6$-tuple $\left(\calG,X_\calG,s_\calG,t_\calG,\iota_\calG,j_\calG\right)$, where $\calG$ and $X_\calG$ are two sets (called respectively the {\em set of arrows} and the {\em set of points or (more commonly) objects}, $s_\calG$ and $t_\calG$ are two surjective maps from $\calG$ to $X_\calG$ (called respectively the {\em source map or source} and the {\em target map or target}), $\iota_\calG$ is a map from $X_\calG$ to $\calG$ (called the {\em unit map}) and $j_\calG$ is a map from $\calG$ to itself (called the {\em inversion}); furthermore, introducing the subset of $\calG\times\calG$ of {\em composable arrows}, denoted by $\calG_2$, defined via
\[
\calG_2=\left\{(g_1,g_2)\in\calG\times\calG\colon s_\calG(g_1)=t_\calG(g_2)\right\},
\] 
there is an operation from $\calG_2$ to $\calG$, the {\em product of the groupoid $\calG$},
\[
\calG_2\ni(g_1,g_2)\mapsto g_1g_2.
\]

The following axioms must be satisfied
\begin{itemize}
\item[i)] for any composable couple $(g_1,g_2)\in\calG_2$, it holds
\[
s_\calG\!(g_1g_2)=s_\calG\!(g_2),\quad t_\calG\!(g_1g_2)=t_\calG\!(g_1).
\] 
\item[ii)] {\bf (Identity axiom)} For any $x\in X_\calG$, it holds
\[
s_\calG(\iota_\calG(x))=x=t_\calG(\iota_\calG(x));
\] 
furthermore, for any $g\in\calG$, it holds
\[
\iota_\calG(t_\calG(g))g=g=g\iota_\calG(s_\calG(g)).
\]
\item[iii)] {\bf (Inversion axiom)} For any $g\in\calG$, it holds
\[
s_\calG(j_\calG(g))=t_\calG(g),\quad t_\calG(j_\calG(g))=s_\calG(g);
\] 
furthermore, the following identities must hold
\[
gj_\calG(g)=\iota_\calG(t_\calG(g)),\quad j_\calG(g)g=\iota_\calG(s_\calG(g)).
\]
\item[iv)] {\bf (Associativity)} For any {\em composable triple $(g_1,g_2,g_3)$}, i.e.\ any triple obeying
\[
t_\calG(g_3)=s_\calG(g_2),\quad t_\calG(g_2)=s_\calG(g_1),
\] 
the identity has to be satisfied
\[
(g_1g_2)g_3=g_1(g_2g_3).
\]
\end{itemize}
\end{Def}

\begin{Rem}
Working in the categorical framework, one could speak of a groupoid as of a category $\calG$, whose morphisms are all invertible.
\end{Rem}

I introduce the following notations: for any two points (or objects) $x,y$ of $X_\calG$, the set $\calG_{x,y}$ is defined via
\[
\calG_{x,y}\colon=\left\{g\in \calG\colon s_\calG(g)=x,\quad t_\calG(g)=y\right\}.
\]
Furthermore, the fibre at $x\in\calG$ of the source map $s_\calG$, resp.\ of the target map $t_\calG$, is denoted by $\calG_{x,\bullet}$, resp.\ $\calG_{\bullet,x}$.
Observe that, for any $x\in X_\calG$, the set $\calG_{x,x}$ is a group, called the {\em isotropy group at $x$}: its multiplication is well-defined, as, for any $x\in X_\calG$, $\calG_{x,x}\subset \calG_2$, it is associative. 
There is also a unit element, which is simply $\iota_\calG(x)$; the inverse of an element $g$ is clearly $j_\calG(g)$.

For the sake of simplicity, a groupoid is denoted simply by $\calG$ (i.e.\ by its set of arrows), instead of writing down the complete $6$-tuple $\left(\calG,X_\calG,s_\calG,t_\calG,\iota_\calG,j_\calG\right)$.

The concept of homomorphism of groupoids (or simply morphism of groupoids) is also needed.
\begin{Def}\label{def-morgroupoid}
Given two groupoids $\calG$ and $\calH$, a {\em homomorphism from $\calG$ to $\calH$} (or simply a morphism from $\calG$ to $\calH$), consists of a couple $\left(\Phi,\varphi\right)$, where $i)$ $\Phi$ is a map from the set of arrows $\calG$ to the set of arrows $\calH$, and $ii)$ $\varphi$ is a map from the set of objects $X_\calG$ to the set of points $X_\calH$, obeying the following requirements:
\begin{itemize}
\item[i)] {\bf (Compatibility between the groupoid structures)} the three diagrams must commute
\begin{equation}\label{eq-commdiagmor}
\begin{CD} 
\calG  @>\Phi>> \calH\\
@Vs_\calG VV              @VV s_\calH V\\
X_\calG @>\varphi>> X_\calH
\end{CD}\quad,\quad 
\begin{CD} 
\calG  @> \Phi >> \calH\\
@V t_\calG VV              @VV t_\calH V\\
X_\calH          @>\varphi>> X_\calH
\end{CD}\quad\text{and}\quad
\begin{CD} 
X_\calG  @> \varphi >> X_\calH\\
@V \iota_\calG VV              @VV \iota_\calH V\\
\calH          @>\Phi>> \calH
\end{CD}\quad.
\end{equation}
\item[ii)] {\bf (Homomorphism property)} For any composable pair $(g_1,g_2)\in\calG_2$, the identity must hold
\begin{equation}\label{eq-homomor}
\Phi(g_1g_2)=\Phi(g_1)\Phi(g_2).
\end{equation} 
\end{itemize}
\end{Def}

\begin{Rem}
In the categorical language, a morphism from the groupoid $\calG$ to the groupoid $\calH$ is a functor between the two categories.
\end{Rem}

\begin{Rem}
Notice that the commutativity of the diagrams (\ref{eq-commdiagmor}) and Identity (\ref{eq-homomor}) imply together that
\[
\Phi\circ j_\calG=j_\calH\circ\Phi.
\]
\end{Rem}
                                
Now, having introduced, the notion of groupoid, we are ready to introduce and discuss the notion of Lie groupoid.
\begin{Def}\label{def-liegroupoid}
A {\em Lie groupoid $\calG$} is a groupoid in the sense of Definition~\ref{def-groupoid}, such that the set of objects $X_{\calG}$ has the structure of a smooth manifold (which has to be Hausdorff as a topological space) and the set of arrows $\calG$ has the structure of a smooth (but perhaps not Hausdorff and even not second-countable) manifold; moreover, the source map $s_{\calG}$ has to be a smooth epimorphism (i.e.\ a surjective map with surjective tangent map at each point), with Hausdorff fibres, and all other structure maps are smooth maps. 
Accordingly, the set of arrows is now called the {\em manifold of arrows}, while the set of objects is called the {\em manifold of objects}.
\end{Def}

\begin{Rem}\label{rem-liegroupoid}
Notice that, for a Lie groupoid $\calG$, the unit map $\iota_{\calG}$ is smooth.
The Identity Axiom for the groupoid $\calG$ implies immediately that also the target map is surjective; moreover, it follows, from the smoothness of $\iota_{\calG}$, that the target map is also a smooth submersion.
Moreover, since both source map and target map are surjective submersions, it follows that the set of composable ``arrows'' $\calG_2$ inherits the structure of a smooth manifold, since it is the restriction to the diagonal of the product $X_{\calG}\times X_{\calG}$ of the smooth manifold
\[
s_{\calG}^{-1}\!\left(X_{\calG}\right)\times t_{\calG}^{-1}\!\left(X_{\calG}\right).
\]
\end{Rem}

According to Definition~\ref{def-liegroupoid}, a homomorphism between Lie groupoids is a homomorphism in the sense of Definition~\ref{def-morgroupoid}, where both maps of the pair $\left(\Phi,\varphi\right)$ are smooth maps.

%We say that a morphism $\left(\Phi,\varphi\right)$ between Lie groupoids $\calG$ and $\calH$ is a {\em submersion}, if $\Phi$ is a submersion from $\calG$ to $\calH$ ; we notice that the compatibility conditions in Definition~\ref{def-morgroupoid} imply that the second map $\varphi$ is a submersion between the manifolds of points $X_{\calG}$ and $X_{\calH}$.

\subsection{Some examples of Lie groupoids}\label{ssec-examples}
Before going further, let me discuss some examples of Lie groupoids.
\begin{itemize}
\item[a)] Any Lie group $G$ is by itself a Lie groupoid; namely, consider the group $G$ as the manifold of arrows and a point $\star$ as the manifold of objects. 
The source map and the target map are thus trivial, since they associate, to any $g$, the point $\star$; the unit map associates to $\star$ the usual identity of $G$ and the inversion map is simply $g\mapsto g^{-1}$, the inverse in group-theoretic sense.
The group axioms ensure that $G$ with the above structure is a groupoid. 

\item[b)] If $G$ is a Lie group and $M$ is a manifold acted on smoothly from the left by $G$, define the {\em action groupoid} as the product $G\times M$ as manifold of arrows and $M$ as the manifold of objects.
The source map is simply the projection onto the second factor, while the target map is given by the left action; finally, the multiplication is defined via the assignment
\[
(g_1,m_1)(g_2,m_2)\colon=(g_1g_2,m_1),\quad\forall g_1,g_2\in G,m_1,m_2\in M.
\]
The action groupoid associated to a group $G$ and a left $G$-set $M$ is commonly denoted by $G\ltimes M$.

\item[c)] The {\em fundamental groupoid $\Pi\!(M)$ over a manifold $M$} is defined as follows: the manifold $M$ itself is the set of points.
For any two points $x,y$, the set of arrows $\Pi!(M)_{x,y}$ from $x$ to $y$ is the set of all homotopy classes of paths from $x$ to $y$, relative to endpoints (thus, the isotropy group at $x$ of the fundamental groupoid of $M$ is the fundamental group $\pi_1\!\left(M,x\right)$ based at the point $x$).
The source map and the target map of the fundamental groupoid are then obvious; its multiplication is in turn induced simply by the composition of composable paths, which is compatible with homotopies fixing endpoints.

\item[d)] Given a manifold $M$, there is a natural Lie groupoid associated to $M$, namely the {\em product groupoid of $M$ with itself}: the manifold of arrows is the product manifold $M\times M$, while the manifold of objects is $M$ itself.
The source map is given by projection onto the second factor, while the target map is given by projection onto the first factor; the unit map is simply the diagonal immersion of $M$ into the product of $M$ with itself.
Multiplication is then naturally given by
\[
(x_1,x_2)(x_2,x_3)\colon=(x_1,x_3).
\] 

\item[e)] The {\em gauge groupoid $\calG(P)$ associated to a principal $G$-bundle $P$ over the manifold $M$} is defined as follows: consider the orbit space of the diagonal action of $G$ on the product $P\times P$ as the manifold of arrows, and the base manifold $M$ of $P$ as the manifold of points.
The target map is given by the composition of the projection onto the first factor with the projection $\pi$ from $P$ to $M$; the source map is in turn induced by the composition of the projection onto the second factor with the map $\pi$.
The product is defined so, that the quotient map from the product groupoid $P\times P$ onto the gauge groupoid is a homomorphism of Lie groupoids; without going into the details, let me just say that the product is explicitly constructed via the {\em division map of $P$}, for which I refer e.g.\ to~\cite{McK} or~\cite{C1} for more details. 
Let me notice at the end that the manifold of arrows of the gauge groupoid may be identified with the total space of the bundle associated to the left action of $G$ on $P$ induced by the right action of $G$ on $P$; this is useful when discussing {\em principal bundles with structure groupoid the gauge groupoid of $P$}.
\end{itemize}

\section{General constructions for Lie groupoids: product groupoid, opposite groupoid and groupoid actions}\label{sec-groupoidconstr}
In this Section, I display some general constructions in the theory of groupoids; in particular, I discuss the concept of product groupoid, and, in more details, the concept of left and right $\calG$-spaces, for a general groupoid $\calG$.
In particular, I introduce the notion of {\em generalized conjugation for groupoids}: it is well-known that it is not possible to define conjugation for general groupoids (due to the fact that not all arrows are composable), but as we will see later that for a general groupoid $\calG$ it is possible nonetheless to define two distinct actions of the product groupoid of $\calG$ with itself on $\calG$, both inducing the usual conjugation on each isotropy group $\calG_x$.
Finally, I introduce the concept of (twisted) equivariant maps between left (and right) groupoid spaces, where the actions may come from distinct groupoids; this is the main notion that I need in order to study equivariant morphisms between principal  bundles with structure groupoid from the point of view of generalized gauge transformations.

Let me end the introduction to the topics of this section with a caveat: 

\fbox{\parbox{12cm}{\bf From now on, every groupoid $\calG$, $\calH$ is meant to be a Lie groupoid; I will explicitly specify if otherwise.}} 

\subsection{The product groupoid of two groupoids $\calG$, $\calH$}\label{ssec-prodgroupoid}
Given two groupoids $\calG$ and $\calH$, with respective source, target, unit maps and inversions, we may form the {\em product groupoid of $\calG$ and $\calH$} by setting
\begin{itemize}
\item[i)] the product manifold $\calG\times \calH$ as the manifold of arrows of the product groupoid;
\item[ii)] the product manifold $X_\calG\times X_\calH$ as the manifold of objects of the product groupoid; 
\item[iii)] the map
\[
s_{\calG\times \calH}\!\left(g,h\right)\colon=\left(s_\calG(g),s_\calH(h)\right),\quad \forall (g,h)\in \calG\times\calH,
\] 
as the source map of the product groupoid;
\item[iv)] the map
\[
t_{\calG\times \calH}\!\left(g,h\right)\colon=\left(t_\calG(g),t_\calH(h)\right),\quad \forall (g,h)\in \calG\times\calH,
\] 
as the target map of the product groupoid; 
\item[v)] the map 
\[
\iota_{\calG\times \calH}\!(x,y)\colon=\left(\iota_\calG(x),\iota_{\calH}(y)\right),\quad \forall (x,y)\in X_{\calG}\times X_\calH,
\] 
as the unit map of the product groupoid;
\item[vi)] the map 
\[
j_{\calG\times \calH}\!(g,h)\colon=\left(j_\calG(g),j_\calH(h)\right),\quad \forall (g,h)\in \calG\times\calH,
\] 
as the inversion of the product groupoid;
\item[vii)] the partial product of the product groupoid is defined by the assignment
\[
(g_1,h_1)(g_2,h_2)\colon=(g_1g_2,h_1h_2),\quad s_{\calG\times\calH}\!\left(g_1,h_1\right)=t_{\calG\times\calH}\!\left(g_2,h_2\right).
\] 
We notice that the definition of product makes sense by the very definition of the source and target map in the product groupoid.
\end{itemize}
It is immediate to check that all axioms of (Lie) groupoid are satisfied: in particular, the product of two Lie groupoids is again a Lie groupoid, as the product of smooth manifolds is again smooth, and the product of smooth maps is again smooth.
Finally, the product of the source maps is clearly surjective, and, by definition of tangent map, it is clearly a submersion, as both its factors are submersions.

\subsection{Left- and right $\calG$ actions for the groupoid $\calG$}\label{ssec-leftactgroupoid}
Given now a groupoid $\calG$ and a smooth manifold $M$, I want to clarify the notion of left $\calG$-action on $M$; the notion of right $\calG$-action is similar, and I mention it briefly.

\begin{Def}\label{def-groupoidaction}
A left $\calG$-action of the groupoid $\calG$ on the (smooth) manifold $M$ consists of a $3$-tuple $\left(M,J_M,\Psi_M\right)$, where $i)$ $J_M$ is a smooth map from $M$ to the manifold of objects $X_\calG$ of the groupoid $\calG$ (called the {\em momentum of the action} or, more briefly, the {\em momentum}, and $ii)$ $\Psi_M$ is a smooth map from $\calG\times_{J_M} M$ to $M$, where
\[
\calG\times_{J_M} M\colon=\left\{(g,m)\in\calG\times M\colon s_\calG(g)=J_M(m)\right\}.
\] 
It is customary to write
\[
\Psi_M\!\left(g,m\right)\colon= gm
\]
(Usually, one speaks also of left $\calG$-action w.r.t.\ the momentum $J_M$.)

Moreover, the following requirements must hold
\begin{itemize}
\item[i)] \[
J_M\!\left(gm\right)=t_\calG(g),\quad \forall (g,m)\in \calG\times_{J_M}M;
\] 
\item[ii)] \[
g_1\left(g_2 m\right)=\left(g_1 g_2\right)m,\quad \forall (g_1,g_2)\in\calG_2,(g_1g_2,m)\in\calG\times_{J_M}M;
\] 
Observe that Condition $\mathnormal{i})$ implies that the previous identity is well-defined.
\item[iii)] \[
\iota_\calG\left(J_M(m)\right)m=m,\quad \forall m\in M.
\] 
\end{itemize}
\end{Def}

\begin{Rem}
Notice that the set $\calG\times_{J_M} M$ is in fact a manifold, as it is the pull-back w.r.t.\ the momentum $J_M$ of the smooth fibration over $X_{\calG}$ defined by the source map.
\end{Rem}

\begin{Rem}
The definition of right $\calG$-action is similar, the only difference being that one has to switch the r{\^o}les of the source and target maps; consequently, the map $\Psi_M$ goes from the product $M\times_{J_M}\calG$ to $M$, and is denoted by
\[
\Psi_M(m,g)\colon=mg.
\]
Equivalently, a right $\calG$-action is a left $\calG^{\op}$-action, and the switch between the two actions is provided by the inversion $j_\calG$. 
\end{Rem}

One says that a left groupoid action of $\calG$ with momentum $J_M$ on a manifold $M$ is {\em free}, if the following condition holds:
\begin{equation*}
g m=m,\quad s_{\calG}(g)=J_M(m)\Rightarrow g=\iota_{\calG}(J_M(m)).
\end{equation*}
This implies that, for any $x\in X_\calG$, any isotropy group $\calG_{x,x}$ operates freely (in the usual sense) on the fibre $J_M^{-1}\!\left(\left\{x\right\}\right)$.
(The concept of freeness for a right $\calG$-action is similar.)
On the other hand, one says that a left $\calG$-action with momentum $J_M$ on a manifold $M$ is {\em transitive}, if the following requirement holds:
\begin{equation*}
\forall m,\widetilde{m}\in M,\exists g\in\calG_{J_M(m),J_M(\widetilde{m})}\quad\text{such that $\widetilde{m}=gm$}.
\end{equation*}
(The definition of transitivity of a right $\calG$-action is similar.)
 
\begin{Rem}
Notice that, if a left (or right) $\calG$-action with momentum $J_M$ on a manifold $M$ is free and transitive, the transitivity condition may be restated as
\[
\forall m,\widetilde{m}\in M,\exists! g\in\calG_{J_M(m),J_M(\widetilde{m})}\quad\text{such that $\widetilde{m}=gm$}.
\]
\end{Rem}

\subsubsection{The generalized conjugation of $\calG$}\label{sssec-genconjgroupoid}
As I have already remarked at the beginning of the Section, a groupoid does not admit a natural notion of conjugation as a usual Lie group.
In fact, in a usual Lie group, the conjugation by an element $g$ of an element $h$ is given by the formula $ghg^{-1}$. 
The natural notion of conjugation for a groupoid would be then to consider conjugation on any isotropy group $\calG_x$, which is clearly a Lie group, and corresponds naturally to the conjugation for a usual Lie group, since, in this case, any isotropy group is equal to the groupoid itself; but this definition is too restrictive.
In fact, one needs a momentum $J_\calG$ from the manifold of arrows $\calG$ to the manifold of objects $X_\calG$ and a left action map $\Psi_\calG$ from $\calG\times_{J_\calG}\calG$ to $\calG$, obeying the three requirements of Definition~\ref{def-groupoidaction}; the left action map, intuitively, has to take the form, whenever it makes sense,
\[
(g,h)\overset{\Psi_\calG}\mapsto ghg^{-1}.
\]  
The conjugation equation requires, by its very definition, that
\[
t_\calG(h)=s_\calG(g),\quad s_\calG(h)=t_\calG(g^{-1})=s_\calG(g)\Rightarrow t_\calG(h)=s_\calG(h).
\]
Thus, the usual conjugation makes sense only on the isotropy groups $\calG_{x,x}$, for any $x\in X_\calG$.
On the other hand, for any Lie group $G$, it is possible to construct four distinct actions of the product $G\times G$ on $G$ itself, namely
\begin{align*}
((g_1,g_2),g_3)&\mapsto g_1g_3g_2^{-1},\\
((g_1,g_2),g_3)&\mapsto g_2g_3g_1^{-1},\\
(g_3,(g_1,g_2))&\mapsto g_1^{-1}g_3g_2,\\
(g_3,(g_1,g_2))&\mapsto g_2^{-1}g_3g_1.
\end{align*}
All four actions are clearly smooth; the two first actions are left actions, while the remaining two are right actions.
There is a natural subgroup of the product $G\times G$, namely the diagonal subgroup, which is naturally isomorphic to $G$; when restricting the two first actions of $G\times G$ on the diagonal subgroup, one gets the same action of $G$ on $G$, which is the left conjugation of $G$, as the following easy computation shows
\[
(g,h)\mapsto ((g,g),h)\mapsto ghg^{-1}=\conj(g)h.
\]
Similarly, the restriction to the diagonal subgroup of the two right actions of $G\times G$ equals the conjugation of $G$ composed with the inversion $g\mapsto g^{-1}$, which is the right conjugation of $G$.
Therefore, the (left or right) conjugation of $G$ can be viewed as particular cases of two more general actions of the product $G\times G$ on $G$, which I call the {\em generalized conjugation of $G$}.

As the following arguments show, the generalized conjugation of groups admits a natural extension to Lie groupoids, which I also call the {\em generalized conjugation of Lie groupoids}.
The first ingredient one needs is a {\em momentum} for the action:
\[
J_{\conj}\!\left(g\right)\colon=\left(t_\calG(g),s_\calG(g)\right),\quad \forall g\in\calG.
\]
Consequently, the manifold $\calG^2\times_{J_{\conj}}\calG$, where the action makes sense, has the form
\[
\calG^2\times_{J_{\conj}}\calG=\left\{\left(g_1,g_2;g_3\right)\in\calG^3\colon \begin{cases}
s_\calG(g_1)&=t_\calG(g_3)\\
s_\calG(g_2)&=s_\calG(g_3)
\end{cases}\right\}.
\]
It makes thus sense to define a map $\Psi_{\conj}$ from $\calG^2\times_{J_{\conj}}\calG$ to $\calG$ as follows:
\begin{equation}\label{eq-genconjgroup}
\Psi_{\conj}\!\left(g_1,g_2;g_3\right)\colon= g_{1}g_3 g_2^{-1},
\end{equation}
where I set for simplicity $g_2^{-1}\colon=j_\calG(g_2)$.

\begin{Prop}\label{prop-genconj}
The triple $\left(\calG^2,J_{\conj},\Psi_{\conj}\right)$ defines a left $\calG^2$-action on $\calG$, which we call the {\em generalized conjugation of $\calG$}. 
\end{Prop}
\begin{proof}
First of all, we notice that the maps $J_{\conj}$ and $\Psi_{\conj}$ are smooth on their domains of definitions.

We then compute, for any triple $\left(g_1,g_2;g_3\right)$ in $\calG^2\times_{J_{\conj}}\calG$, the following expression:
\begin{align*}
\left(J_{\conj}\circ\Psi_{\conj}\right)\!\left(g_1,g_2;g_3\right)&=J_{\conj}\!\left(g_1g_3g_2^{-1}\right)=\\
&\overset{\text{by definition of $J_{\conj}$}}=\left(t_{\calG}\!\left(g_1g_3g_2^{-1}\right),s_{\calG}\!\left(g_1g_3g_2^{-1}\right)\right)=\\
&=\left(t_{\calG}(g_1),s_{\calG}(g_2^{-1})\right)=\\
&=\left(t_{\calG}(g_1),t_{\calG}(g_2)\right)=\\
&=t_{\calG^2}\!\left(g_1,g_2\right), 
\end{align*}
which proves the first requirement for $\left(\calG^2,J_{\conj},\Psi_{\conj}\right)$ to be a left $\calG^2$ action.

Second, we compute explicitly
\begin{align*}
\Psi_{\conj}\!\left(g_1,g_2;\Psi_{\conj}\!(\left(h_1,h_2;g_3\right)\right)&=\Psi_{\conj}\!\left(g_1,g_2;h_1g_3h_2^{-1}\right)=\\
&=g_1\left(h_1g_3h_2^{-1}\right)g_2^{-1}=\\
&=(g_1h_1)g_3\left(g_2h_2\right)^{-1}=\\
&=\Psi_{\conj}\!\left(g_1h_1,g_2h_2;g_3\right),
\end{align*}
whenever the identity makes sense.

Finally, we compute
\begin{align*}
\Psi_{\conj}\!\left(\iota_{\calG^2}\!\left(J_{\conj}(g)\right);g\right)&=\Psi_{\conj}\!\left(\iota_\calG(t_\calG(g)),\iota_\calG(s_\calG(g));g\right)=\\
&=\iota_\calG(t_\calG(g))g\iota_\calG(s_\calG(g))^{-1}=\\
&=g,\quad \forall g\in \calG,
\end{align*}
which ends the proof of the Proposition.
\end{proof}

\begin{Rem}\label{rem-genconj}
There is a similar (still distinct!) left $\calG^2$-action on $\calG$; in fact, one could consider the map $\overline{J}_{\conj}$ from $\calG$ to $X_\calG\times X_\calG$ given by
\[
\overline{J}_{\conj}(g)\colon=\left(s_{\calG}(g),t_{\calG}(g)\right),
\]
whence
\[
\calG^2\times_{\overline{J}_{\conj}}\calG=\left\{(g_1,g_2;g_3)\in\calG^3\colon \begin{cases}
s_\calG(g_1)&=s_{\calG}(g_3)\\
s_{\calG}(g_2)&=t_\calG(g_3)
\end{cases}
\right\},
\]
and the map $\overline{\Psi}_{\conj}$ from $\calG^2\times_{\overline{J}_{\conj}}\calG$ to $\calG$ via
\[
\overline{\Psi}_{\conj}\!\left(g_1,g_2;g_3\right)\colon=g_2g_3g_1^{-1}.
\]
It is not difficult to verify that the triple $\left(\calG,\overline{J}_{\conj},\overline{\Psi}_{\conj}\right)$ defines also a left $\calG^2$-action on $\calG$.
\end{Rem}

\begin{Rem}\label{rem-conjright}
Notice that the maps $J_{\conj}$ and $\overline{J}_{\conj}$ define also right $\calG^2$-actions on $\calG$, the {\em right generalized conjugations}: namely, on the set $\calG\times_{J_{\conj}}\calG^2$, resp.\ $\calG\times_{\overline{J}_{\conj}}\calG^2$, we define the map $\Psi_{\conj}^R$, resp.\ $\overline{\Psi}_{\conj}^R$, by the formula
\begin{align*}
(g_3;g_1,g_2)&\overset{\Psi_{\conj}^R}\mapsto g_1^{-1}g_3g_2,\quad\text{resp.}\\
(g_3;g_1,g_2)&\overset{\overline{\Psi}_{\conj}^R}\mapsto g_2^{-1}g_3g_1.
\end{align*}
\end{Rem}

\subsection{Twisted equivariant maps between groupoid-spaces}\label{ssec-equivgroupoid}
I define and discuss briefly the concept of {\em equivariant map between groupoid-spaces}. 
For simplicity, by groupoid-space, I mean here a manifold $M$ acted on from the left by a groupoid $\calG$.

For our purposes, I will consider the most general situation, namely a left $\calG$-space $\left(M,J_M,\Psi_M\right)$ and a left $\calH$-space $\left(N,J_N,\Psi_N\right)$, where $\calG$, $\calH$ are two groupoids and $M$, $N$ are two manifolds
.
\begin{Def}\label{def-equivgroupoid}
A {\em (twisted) equivariant map} between the left $\calG$-space $M$ and the left $\calH$-space $N$ consists of a triple $\left(\Theta,\Phi,\varphi\right)$, where $\Theta$ is a smooth map from the manifold $M$ to the manifold $N$, and the pair $\left(\Phi,\varphi\right)$ is a morphism from the groupoid $\calG$ to the groupoid $\calH$ in the sense of Definition~\ref{def-morgroupoid}.

Moreover, the following two diagrams must commute:
\begin{itemize}
\item[i)] \[
\begin{CD} 
M  @>\Theta>> N\\
@VJ_M VV              @VV J_N V\\
X_\calG @>\varphi>> X_\calH
\end{CD}\quad;
\] 
\item[ii)] \[ 
\begin{CD}
\calG\times_{J_M}M  @>\Phi\times \Theta  >> \calH\times_{J_N}N\\
@V\Psi_M VV              @VV \Psi_N V\\
M @>\Theta>> N 
\end{CD}
\]
\end{itemize}
\end{Def}  

\begin{Rem}
The first commutative diagram in Definition~\ref{def-equivgroupoid} implies that $\Phi\times \Theta$ maps really the manifold $\calG\times_{J_M}M$, where the $\calG$-action is well-defined, to the manifold $\calH\times_{J_N}N$, where the $\calH$-action is well-defined, as the following explicit computation shows:
\begin{align*}
J_N\!\left(\Theta(m)\right)&=\varphi(J_M(m))=\\
&=\varphi(s_\calG(g))=\\
&=s_{\calH}(\Phi(g)),\quad \forall (g,m)\in \calG\times_{J_M} M.
\end{align*} 

Usually, the second diagram may be rewritten as the identity:
\[
\Theta(gm)=\Phi(g)\Theta(m),\quad \forall (g,m)\in \calG\times_{J_M} M,
\]
which corresponds clearly to the usual definition of (twisted by $\Phi$) equivariance of a map $\Theta$ from a left $G$-space to a left $H$-space, for $G$, $H$ usual groups.
\end{Rem}

The concept of (twisted) equivariant map between right groupoid-spaces is similar, the only difference being that one has to invert the factors in the products $\calG\times_{J_M}M$, $\calH\times_{J_N}N$ and $\Phi\times \Theta$.

\section{Principal bundles with structure groupoid}\label{sec-princbun}
An important notion in differential geometry is that of {\em principal bundle}: a principal bundle $P$ with structure group $G$ over the manifold $M$ is a triple $\left(P,\pi,M\right)$, where $P$ and $M$ are both smooth manifolds, $\pi$ is a surjective submersion from $P$ to $M$ (i.e.\ a map whose tangent map at any point $p$ of $P$ is surjective) such that the following requirements hold:
\begin{itemize}
\item[i)] the group $G$ acts freely from the right on $P$; 
\item[ii)] the projection $\pi$ is $G$-invariant:
\[
\pi(pg)=\pi(p),\quad \forall p\in P,g\in G;
\] 
\item[iii)] $P$ is {\em locally trivial} in the following sense: given a point $x\in M$, there exists an open neighbourhood $U=U_x$ of $x$ in $M$ and a diffeomorphism $\varphi_U$ 
\[
\varphi_U\colon \pi^{-1}\!\left(U\right)\mapsto U\times G,
\]
and $\varphi_U$ is $G$-equivariant
\[
\phi_U(pg)=\phi_U(p)g,
\]
where $G$ acts from the right on $U\times G$ by right multiplication on the second factor of any pair in $U\times G$, and satisfies the equation
\[
\pr_1\circ \varphi_U=\pi,
\]
where $\pr_1$ denotes projection onto the first factor of any pair in $U\times G$.
\end{itemize}
The trivial principal bundle over $M$ is simply the triple $(M\times G,\pr_1,M)$, where $G$ acts from the right on the product manifold $M\times G$ by right multiplication on the second factor of any pair.

I give now the notion of principal bundle with groupoid structure, namely, I want to define an analogue of principal bundles in the above sense, where I replace the structure group $G$ by a more general groupoid $\calG$.
The natural concepts appearing in the theory of usual principal bundles that may be translated immediately to the theory of principal bundles with structure groupoid are that of right $\calG$-space and of surjective submersion; it remains therefore to give a criterion which in some sense mimics the ``triviality condition''.

\begin{Def}\label{def-princgroupoid}
A {\em principal bundle $P$ with groupoid structure $\calG$ over the manifold $M$} is a $4$-tuple $\left(P,\pi,\varepsilon,M\right)$, where $i)$ $P$ and $M$ are smooth manifolds and $ii)$ the pair $(P,\varepsilon)$ defines a structure of right $\calG$-space on $P$ (we drop the right action map, denoting it simply by a product or, when needed, by $\Psi$).

Moreover, the following requirements must hold:
\begin{itemize}
\item[i)] the map $\pi$ is a surjective submersion from $P$ to $M$. 
\item[ii)] The map $\pi$ is $\calG$-invariant, i.e.\ the following diagram commutes
\[
\begin{CD} 
P\times_{\varepsilon}\calG  @>\Psi>> P \\
@V\pr_1 VV  @VV\pi V\\
P @>\pi>> M
\end{CD}\quad .
\] 
\item[iii)] The map $\left(\pr_1,\Psi\right)$ defined via
\begin{align*}
\left(\pr_1,\Psi\right)\colon P\times_\varepsilon \calG&\to P\times_M P,\\ 
(p,g)&\mapsto (p,pg),
\end{align*}
is a diffeomorphism; by $P\times_M P$, we mean
\[
P\times_M P\colon=\left\{(p,q)\in P\times P\colon \pi(p)=\pi(q)\right\}.
\]
\end{itemize} 
\end{Def}

\begin{Rem}
The notion of principal bundles with structure groupoid as in the previous definition is not new: in fact, it was introduced by Connes in~\cite{Con} for studying the holonomy groupoid of a foliation, and used extensively later by H{\ae}fliger in~\cite{Haef}, although they used a local description in terms of nonabelian Cech cohomology for groupoids (still, their notion of nonabelian Cech cohomology for groupoids, although correct, lacks of an explicit mentioning of what I call local momenta; I plan to return to this point in subsequent works).
Later, M{\oe}rdijk~\cite{Moer1} took a different point of view, working nonlocally, introducing the notion of {\em cocycle on $M$ with values in $\calG$}, which mentions explicitly the presence of a momentum, and which corresponds, in terms of groupoids, to the division map for ordinary principal bundles discussed extensively by MacKenzie~\cite{McK}.
Finally, the nonlocal point of view was formulated in a definitive way by Mrcun~\cite{Mrcun} and M{\oe}rdijk~\cite{Moer2}, which is the point of view that I take here.
Let me notice that the local point of view, in terms of nonabelian Cech cohomology, has also many advantages, among them, e.g., the possibility of constructing explicitly many examples of principal bundles; still, I will only mention briefly the local nature of principal bundles with structure groupoids here, devoting subsequent works to this aspect of the theory. 
\end{Rem}

\begin{Exa}
(For more examples and details about them, I refer to~\cite{Moer2},~\cite{C2})
\begin{itemize}
\item[i)] Given a Lie group $G$, considered as a groupoid over a point with trivial target, source and unit map, and a manifold $M$, a principal bundle $P$ with structure groupoid $G$ is the same as a principal $G$-bundle in the usual sense. 
\item[ii)] Given a Lie group $G$ acting from the left on a manifold $M$, one can consider the action groupoid $G\ltimes M$.
Then, principal $G\ltimes M$-bundles $P$ over $X$ are in one-to-one correspondence with principal $G$-bundles over $X$ in the usual sense with a global section of the associated bundle $P\times_G M$ over $X$. 
\item[iii)] Consider a manifold $M$ and the product groupoid $\Pi(M)$, and another manifold $X$.
Then, there is a unique principal $\Pi(M)$-bundle over $X$, namely the $4$-tuple $\left(X\times M,\pr_X,\pr_M,X\right)$.  
\end{itemize}
\end{Exa}

\begin{Rem}\label{rem-freeness}
The meaning of the third axiom is that the groupoid $\calG$ operates {\em freely and transitively on each fibre of $\pi$}.
In fact, assume the identity holds
\[
pg=p,\quad p\in P,g\in\calG\quad\text{such that $t_{\calG}(g)=\varepsilon(p)$}.
\]
It follows that both pairs $(p,g)$ and $(p,\iota_{\calG}(\varepsilon(p)))$, both in $P\times_\varepsilon \calG$, are mapped by the diffeomorphism $(\pr_1,\Psi)$ to the same image, namely $(p,p)$; hence, 
\[
g=\iota_{\calG}(\varepsilon(p)).
\] 

If, on the other hand, we take any two points $p$ and $q$ of $P$, lying in the same fibre of $\pi$, we have $(p,q)\in P\times_M P$;
since $\left(\pr_1,\Psi\right)$ is a diffeomorphism, it follows that 
\[
q=pg,\quad g\in\calG\quad \varepsilon(q)=s_{\calG}(g),\quad \varepsilon(p)=t_{\calG}(g),
\]
whence $g\in \calG_{\varepsilon(q),\varepsilon(p)}$.
\end{Rem}

The smooth inverse of the map $\left(\pr_1,\Psi\right)$ is usually denoted by $\Phi_P$.
It follows by its very definition:
\[
P\times_M P\ni (p,q)\overset{\Phi_P}\mapsto(\Phi_{P,1}(p,q),\Phi_{P,2}(p,q))\in P\times_\varepsilon \calG,
\]
whence $\Phi_{P,1}$, resp.\ $\Phi_{P,2}$, is a smooth map from $P\times_M P$ to $P$, resp.\ $\calG$; clearly, 
\[
\Phi_{P,1}(p,q)=p,\quad\forall (p,q)\in P\times_M P.
\]
On the other hand, since the image of $\Phi_P$ lies in $P\times_\varepsilon \calG$, it follows
\begin{align*}
\varepsilon\left(\Phi_{P,1}(p,q)\right)&=\varepsilon(p)=t_{\calG}(\Phi_{P,2}(p,q)).
\end{align*}
I denote, from now on, the map $\Phi_{P,2}$ simply by $\phi_P$.
I will analyze its properties in detail later; notice that the function $\phi_P$ was already introduced by M{\oe}rdijk in~\cite{Moer1}, where the pair $\left(\phi_P,\varepsilon\right)$ was called a cocycle on $M$ with values in $\calG$, and earlier by MacKenzie in~\cite{McK} in the case of ordinary principal bundles with structure group $G$, where it was called the {\em division map of $P$}.
I prefer the notation $\phi_P$ in order to make explicit its dependence on the principal bundle $P$.

\subsection{The unit bundle of a groupoid $\calG$ and the trivial bundle}\label{ssec-unitbun}
In this Subsection, I consider two important examples of principal $\calG$-bundles, namely the {\em unit bundle of $\calG$} and the {\em trivial principal $\calG$-bundle}.

As the readers have surely noticed, there is no trace (apparently) in Definition~\ref{def-princgroupoid} of the triviality condition present in the definition of a principal $G$-bundle. 
This is because, in fact, the definition of trivial principal $\calG$-bundle requires more care that the definition of the usual trivial $G$-bundle, and requires also the notion of unit bundle; nonetheless, we will see later that some sort of triviality condition holds also for principal $\calG$-bundles.

\begin{Rem}
I will later discuss more carefully the ``local triviality problem'' for principal bundles with structure groupoids: I will namely prove an equivalence between Definition~\ref{def-princgroupoid} and local data obeying some cochain properties.
In fact, the second characterization provides a useful way for constructing non-trivial principal bundles with groupoid structure.
\end{Rem}

\begin{Def}\label{def-unitbun}
The {\em unit bundle of the Lie groupoid $\calG$} consists of the $4$-tuple $\left(\calG,t_\calG,s_\calG,X_\calG\right)$ (thus, it is a bundle over the manifolds of objects of $\calG$), and the right $\calG$-action on itself is given by right multiplication; it is usually denoted by $\mathcal{U}_{\calG}$.
\end{Def}
Notice that the choice of right multiplication on $\calG$ as right $\calG$-action on the unit bundle is in accordance with our choice of the map $\varepsilon$ in the previous definition.
Notice also that the unit bundle $\calU_\calG$ is a principal $\calG$-bundle in the sense of Definition~\ref{def-princgroupoid}: namely, since $\calG$ is a Lie groupoid, the target map is also a smooth submersion, whence $i)$ of Definition~\ref{def-princgroupoid} is satisfied.

It remains to show that also $ii)$ holds; I show this explicitly to motivate the terminology ``division map''.
The map $\left(\pr_1,\Psi\right)$, where $\Psi$ is the right multiplication map, takes the explicit form
\[
\calG\times_{s_\calG}\calG\ni (g,h)\mapsto (g,gh)\in\calG\times_{t_{\calG}}\calG.
\]
It is then easy to prove that the previous map has an explicit smooth inverse, which turns out to be
\[
\calG\times_{t_\calG}\calG\ni (g,h)\mapsto (g,g^{-1}h)\in \calG\times_{s_{\calG}}\calG;
\] 
hence, the division map $\phi_{\calU_{\calG}}$ associated to the unit bundle is simply given by
\[
\phi_{\calU_{\calG}}(g,h)=g^{-1}h,
\] 
which is the nonabelian version for groupoids of the usual division map for abelian groups.
The unit bundle $\calU_\calG$ is also called the {\em trivial $\calG$-bundle over $X_\calG$}.

%\begin{Rem}
%The example of the unit bundle $\calU_\calG$ clarifies the name {\em division map} chosen by MacKenzie.
%\end{Rem}

In order to give the definition of trivial principal $\calG$-bundle over a manifold $M$, one needs the notion of pull-back bundle.
\begin{Def}\label{def-pullbun}
If the $4$-tuple $\left(P,\pi,\varepsilon,N\right)$ is a principal $\calG$-bundle over $N$ and $M\overset{f}\to N$ is a smooth map from the manifold $M$ to the manifold $N$, the {\em pull-back $f^*P$ of $P$ w.r.t.\ $f$} is defined via
\[
f^*P\colon=\left\{(m,p)\in M\times P\colon f(m)=\pi(p)\right\}.
\]
\end{Def}
Considering the $4$-tuple $\left(f^*P,\pr_1,\varepsilon\circ\pr_2,M\right)$, then one can prove that it defines a principal $\calG$-bundle over $M$, where $\pr_i$, $i=1,2$, denotes projection onto the $i$-th term of $f^*P$.
In fact, it is easy to verify that $\pr_1$ is a surjective submersion; the right $\calG$ action is defined along the map
\[
(m,p)\overset{\varepsilon\circ\pr_2}\to \varepsilon(p)
\]
and takes the explicit form
\[
f^*P\times_{\varepsilon\circ\pr_2}\calG\ni (m,p;g)\mapsto(m,pg).
\]
If two points $(m_1,p)$ and $(m_2,q)$ of $f^*P$ belong to the same fibre, it follows
\[
m_1=m_2\Rightarrow f(m_1)=\pi(p)=f(m_2)=\pi(q)\Leftrightarrow q=p\phi_P(p,q). 
\] 
Thus, the map 
\[
f^*P\times_{\varepsilon\circ\pr_2}\calG\ni (m,p;g)\to (m,p;m,pg)
\]
is a diffeomorphism, where the smooth inverse is given explicitly by
\[
f^*P\times_{M} f^*P\ni(m,p;m,q)\to \left(m,p;\phi_P(p,q)\right),
\]
using explicitly the division map $\phi_P$ of $P$.

\begin{Def}\label{def-trivbun}
Given a groupoid $\calG$ and a smooth map $\alpha$ from a manifold $M$ to the manifold of objects $X_{\calG}$ of the groupoid $\calG$, one may consider the pull-back bundle $\alpha^*\calU_{\calG}$ of the unit bundle of $\calG$.
By its very definition, the total space of this bundle has the form
\[
\alpha^*\calU_{\calG}=\left\{(m,g)\in M\times \calG\colon \alpha(m)=t_{\calG}(g)\right\}.
\]
The bundle $\alpha^*\calU_{\calG}$ is called the {\em trivial $\calG$-bundle over $M$ w.r.t.\ $\alpha$}. 
\end{Def}

\begin{Rem}\label{rem-trivial}
Notice that, while there is only one trivial principal $G$-bundle over a manifold $M$, with $G$ a group, there can be in principle many {\em distinct} trivial $\calG$-bundles over the same base. 
\end{Rem}

\begin{Rem}
Observe that the momentum map $\varepsilon$, along which the right action of $\calG$ on $P$ is defined, is a surjective submersion in the case of a trivial bundle, as it is the composition of two surjective submersions.
\end{Rem}

It is possible to prove a local triviality condition for principal $\calG$-bundles $\left(P,\pi,\varepsilon,M\right)$.
\begin{Prop}
Any principal $\calG$-bundle $\left(P,\pi,\varepsilon,M\right)$ is locally diffeomorphic to a trivial bundle.
\end{Prop} 
\begin{proof}
consider a point $m$ and we choose a local section $\sigma$ of $\pi$ (it is possible, since $\pi$ is a surjective submersion) over an open neighbourhood $U=U_m$, and consider the (smooth) composite map 
\[
\alpha=\alpha_U=\colon=\varepsilon\circ \sigma.
\]
Consider then the map
\[
\alpha^*\calU_{\calG}\ni(m,g)\overset{\varphi_U}\mapsto \sigma(m)g;
\]
by the very definition of the map $\alpha$ and Definition~\ref{def-pullbun}, the map $\varphi_U$ is well-defined and smooth.
Moreover, it has a smooth inverse, which is given by
\[
\pi^{-1}\!(U)\ni p\overset{\psi_U}\mapsto \left(\pi(p),\phi_P\!\left(\sigma(\pi(p)),p\right)\right).
\]
It is clear by the very definition of $\alpha$ that the previous map maps the restriction $\pi^{-1}\!(U)$ of $P$ to $U$ to the trivial bundle $\alpha^*\calU_{\calG}$ and that it is smooth, as a composition of smooth maps.

Let me prove that the map $\psi_U$ is the inverse of $\varphi_U$.
Namely, one has
\begin{align*}
\phi_U\!\left(\psi_U(p)\right)&=\phi_U\!\left(\pi(p),\phi_P\!\left(\sigma(\pi(p)),p\right)\right)=\\
&=\sigma\!(\pi(p))\phi_P\!\left(\sigma(\pi(p)),p\right)=\\
&=p, 
\end{align*}
by Proposition~\ref{prop-prophi}, which will be proved later in Subsection~\ref{ssec-propphi}.
On the other hand, one has
\begin{align*}
\psi_U\!\left(\phi_U(m,g)\right)&=\psi_U\!\left(\sigma(m)g\right)=\\
&=\left(\pi(\sigma(m)g),\phi_P\!\left(\sigma(\pi(\sigma(m)g)),\sigma(m)g\right)\right)=\\
&=\left(m,\phi_P\!\left(\sigma(m),\sigma(m)\right)g\right)=\\
&=(m,g),
\end{align*}
where was used the fact that $\sigma$ is a section of $\pi$, that $\pi$ is $\calG$-invariant and again of Proposition~\ref{prop-prophi}.
\end{proof}

\subsection{Product bundle and fibred product of bundles}\label{ssec-fibprodbun}
In this subsection, I discuss the notion of product of two principal bundles with structure groupoid in the sense of Definition~\ref{def-princgroupoid}: I consider first the case of two principal bundles $P_1$ and $P_2$, over base spaces $M_1$ and $M_2$ and with structure groupoids $\calG_1$ and $\calG_2$ respectively, and I show that there is a natural notion of product $P_1\times P_2$, which can be shown to be a principal bundle over the product of the bases $M_1$ and $M_2$, whose structure groupoid is the product groupoid $\calG_1\times\calG_2$.
I then specialize on the particular case, where the base spaces and the structure groupoids coincide; in this case, it is possible to consider the restriction of the product $P_1\times P_2$ to the diagonal of $M\times M$ ($M$ being the common base space of both $P_1$ and $P_2$, while $\calG$ is the common structure groupoid), and the result is a principal bundle over $M$ with structure groupoid $\calG^2$, the so-called {\em fibred product bundle}, which will play a central r{\^o}le in subsequent computations.

Let me consider first two principal bundles $\left(P_1,\pi_1,\varepsilon_1,M_1\right)$ and $\left(P_2,\pi_2,\varepsilon_2,M_2\right)$, whose structure groupoids are $\calG_1$ and $\calG_2$ respectively.
I consider the product manifold $P_1\times P_2$; there are two natural maps from $P_1\times P_2$ to the product manifolds $M_1\times M_2$ and $X_{\calG_1}\times X_{\calG_2}$ respectively, namely
\begin{equation}\label{eq-prodprojmom}
\pi_1\times \pi_2\colon\begin{cases}
P_1\times P_2&\to M_1\times M_2\\
(p_1,p_2)&\mapsto \left(\pi_1(p_1),\pi_2(p_2)\right)
\end{cases},\quad \varepsilon_1\times \varepsilon_2\colon\begin{cases}
P_1\times P_2&\to X_{\calG_1}\times X_{\calG_2}\\ 
(p_1,p_2)&\mapsto \left(\varepsilon_1(p_1),\varepsilon_2(p_2)\right).
\end{cases}
\end{equation}
Notice that both maps are smooth, as they are products of smooth maps.
\begin{Lem}\label{lem-prodbun}
Given two principal bundles $\left(P_1,\pi_1,\varepsilon_1,M_1\right)$ and $\left(P_2,\pi_2,\varepsilon_2,M_2\right)$ with structure groupoids $\calG_1$ and $\calG_2$ respectively, in the sense of Definition~\ref{def-princgroupoid}, the $4$-tuple $\left(P_1\times P_2,\pi_1\times\pi_2,\varepsilon_1\times\varepsilon_2,M_1\times M_2\right)$ is a principal $\calG_1\times\calG_2$-bundle over $M_1\times M_2$.
\end{Lem}
\begin{proof}
The right $\calG_1\times\calG_2$-action on $P_1\times P_2$ is defined via the momentum $\varepsilon_1\times\varepsilon_2$ as follows: the manifold, where the action is well-defined, is
\[
\left(P_1\times P_2\right)\times_{\varepsilon_1\times\varepsilon_2}\left(\calG_1\times\calG_2\right)\colon=\left\{(p_1,p_2;g_1,g_2)\in P_1\times P_2\times\calG_1\times\calG_2\colon \begin{cases}
\varepsilon_1(p_1)&=t_{\calG_1}(g_1),\\
\varepsilon_2(p_2)&=t_{\calG_2}(g_2)
\end{cases}
\right\},
\]
and the right action map is simply
\[
\left(P_1\times P_2\right)\times_{\varepsilon_1\times\varepsilon_2}\left(\calG_1\times\calG_2\right)\ni (p_1,p_2;g_1,g_2)\mapsto (p_1g_1,p_2g_2)\in P_1\times P_2.
\]
It is immediate to verify that the above rule defines a right $\calG_1\times\calG_2$-action on $P_1\times P_2$.

The bundle projection is, by its very definition, the product of both bundle projections $\pi_1$ and $\pi_2$, hence it is clearly a smooth surjective submersion.
By definition and by the second requirement of Definition~\ref{def-princgroupoid}, it follows also immediately that the product of the bundle projections $\pi_1$ and $\pi_2$ is $\calG_1\times\calG_2$-invariant.

It remains to show the third requirement of Definition~\ref{def-princgroupoid}.
We have to show that the map 
\begin{equation}\label{eq-prodtrans}
\begin{aligned}
&\left(P_1\times P_2\right)\times_{\varepsilon_1\times\varepsilon_2}\left(\calG_1\times\calG_2\right)\ni (p_1,p_2;g_1,g_2)\mapsto\\
&\mapsto(p_1,p_2;p_1g_1,p_2g_2)\in (P_1\times P_2)\times_{M_1\times M_2}(P_1\times P_2)
\end{aligned}
\end{equation}
is a diffeomorphism.
It is a smooth map, as one may view it as the composite map 
\begin{align*}
&\left(P_1\times P_2\right)\times_{\varepsilon_1\times\varepsilon_2}\left(\calG_1\times\calG_2\right)\ni(p_1,p_2;g_1,g_2)\mapsto\\
&\mapsto ((p_1,g_1),(p_2,g_2))\in (P_1\times_{\varepsilon_1}\calG_1)\times(P_2\times{\varepsilon_2}\calG_2)\\
&\mapsto ((p_1,p_1g_1),(p_2,p_2g_2))\in (P_1\times_{M_1} P_1)\times(P_2\times_{M_2}P_2)\\
&\mapsto (p_1,p_2;p_1g_1,p_2g_2)\in (P_1\times P_2)\times_{M_1\times M_2}(P_1\times P_2), 
\end{align*}
and all maps are clearly smooth.
Its inverse is given by the composite map
\begin{align*}
&(P_1\times P_2)\times_{M_1\times M_2}(P_1\times P_2)\ni (p_1,p_2;\tildep_1,\widetilde{p}_2)\mapsto\\ 
&\mapsto((p_1,\tildep_1),(p_2,\widetilde{p}_2))\in (P_1\times_{M_1}P_1)\times \left(P_2\times_{M_2}P_2\right)\\
&\mapsto ((p_1,\phi_{P_1}(p_1,\tildep_1));(p_2,\phi_{P_2}(p_2,\widetilde{p}_2))\in(P_1\times_{\varepsilon_1} \calG_1)\times\left(P_2\times_{\varepsilon_2}\calG_2\right)\\
&\mapsto(p_1,p_2;\phi_{P_1}(p_1,\tildep_1),\phi_{P_2}(p_2,\widetilde{p}_2))\in\left(P_1\times P_2\right)\times_{\varepsilon_1\times\varepsilon_2}\left(\calG_1\times\calG_2\right).
\end{align*}
It is clear that all maps are smooth, thus the map in Equation (\ref{eq-prodtrans}) defines a diffeomorphism, hence completing the proof. 
\end{proof}
I consider now two principal bundles $P_1=P$ and $P_2=\widetilde{P}$, with the same structure groupoid $\calG$ and over the same manifold $M$, whose right $\calG$-actions are defined along the maps $\varepsilon$ and $\widetilde{\varepsilon}$ respectively, and whose projections are $\pi$ and $\widetilde{\pi}$ respectively.
\begin{Cor}\label{cor-prodbun}
The $4$-tuple $\left(P\times\widetilde{P},\pi\times\widetilde{\pi},\varepsilon\times\widetilde{\varepsilon},M\times M\right)$ is a principal bundle over $M\times M$ with structure groupoid $\calG\times \calG$.
\end{Cor}
Now, let me consider the restriction to the diagonal $\Delta_M\subset M\times M$ of the principal bundle $P\times \widetilde{P}$, for any two principal bundles $P$ and $\widetilde{P}$ over the same base $M$ and with the same structure groupoid $\calG$, as in Lemma~\ref{lem-prodbun}; I prove that it is a principal bundle over $M$ with structure groupoid $\calG^2$, which I call the {\em fibred product bundle of $P$ and $\widetilde{P}$}.
\begin{Lem}\label{lem-fibprod}
The $4$-tuple $\left(P\odot\widetilde{P},\overline{\pi},\varepsilon\times\widetilde{\varepsilon},M\right)$, where the manifold $P\odot\widetilde{P}$ is defined by
\[
P\odot\widetilde{P}\colon=\left\{(p,\tildep)\in P\times \widetilde{P}\colon \pi(p)=\widetilde{\pi}(\tildep)\right\},
\]
and the projection $\overline{\pi}$ is 
\[
\overline{\pi}(p,\tildep)=\pi(p)=\widetilde{\pi}(\tildep),
\]
defines a principal $\calG^2$\-bundle over $M$, which is called the fibred product bundle of $P$ and $\widetilde{P}$.
\end{Lem}
\begin{proof}
It is clear that the right action of $\calG^2$ on the product bundle $P\times\widetilde{P}$ restricts to a right $\calG^2$-action on the total space $P\odot\widetilde{P}$ on the fibred product bundle along the same map $\varepsilon\times\widetilde{\varepsilon}$.
It remains to show that the bundle projection $\overline{\pi}$ is a surjective submersion and that the map
\[
(P\odot \widetilde{P})\times_{\varepsilon\times\widetilde{\varepsilon}}(\calG^2)\ni (p,\tildep;g_1,g_2)\overset{(\pr_1^{P\odot\widetilde{P}},\Psi_{P\odot\widetilde{P}})}\mapsto (p,\tildep;pg_1,\tildep g_2)\in (P\odot\widetilde{P})\times_M(P\odot\widetilde{P})
\]
is a diffeomorphism; it is clear that this map is well-defined, as both projections $\pi$ and $\widetilde{\pi}$ are $\calG$-invariant.
The bundle projection $\overline{\pi}$ is obviously smooth and surjective; by its very definition and by the definition of tangent map, it follows also that $\overline{\pi}$ is a submersion.

Finally, it is clear that the map
\[
(P\odot\widetilde{P})\times_M(P\odot\widetilde{P})\ni (p,\tildep;q,\widetilde{q})\mapsto (p,\tildep;\phi_P(p,q),\phi_{\widetilde{P}}(\tildep,\widetilde{q}))\in \left(P\odot \widetilde{P}\right)\times_{\varepsilon\times \widetilde{\varepsilon}}\left(\calG^2\right)
\]
is well-defined and smooth; it is also immediate to check that it is the inverse of the map $\left(\pr_1^{P\odot\widetilde{P}},\Psi_{P\odot \widetilde{P}}\right)$, which is thus a diffeomorphism.
\end{proof}
\begin{Rem}
It is immediate to verify that the fibred product bundle $P\odot\widetilde{P}$ of $P$ and $\widetilde{P}$ can be identified with the pull-back of the product bundle $P\times\widetilde{P}$ w.r.t.\ the diagonal map $\Delta$; nonetheless, I preferred to give a direct proof of all the requirements in Definition~\ref{def-princgroupoid}.
\end{Rem}

\subsection{Properties of the division map $\phi_P$}\label{ssec-propphi}
In this Subsection I analyze in detail the previously introduced division map $\phi_P$, the second component of the inverse of the diffeomorphism $\left(\pr_1,\Psi\right)$.
As already seen, the map $\phi_P$ is defined on $P\times_M P$ and takes its values in the structure groupoid $\calG$ of $P$; it is obvious that one can identify $P\times_M P$ with the total space of the fibred product bundle $P\odot P$, which is, by Lemma~\ref{lem-fibprod} of Subsection~\ref{ssec-fibprodbun}, a right $\calG^2$-space.

Moreover, the map $\phi_P$ satisfies by its very definition the equation
\[
t_{\calG}(\phi_P(p,q))=\varepsilon(p),\forall (p,q)\in P\odot P.
\]

Since $\Phi_P$ is the inverse of $\left(\pr_1,\Psi\right)$, it follows immediately
\begin{align*}
(p,q)&\overset{\Phi_P}\mapsto \left(p,\phi_P(p,q)\right)\\
&\overset{\left(\pr_1,\Psi\right)}\mapsto \left(p,p \phi_P(p,q)\right)\overset{!}=(p,q),\quad \forall(p,q)\in P\odot P,
\end{align*}
whence it follows that the division map $\phi_P$ is defined uniquely by the equation
\begin{equation}\label{eq-defphi}
q=p\phi_P(p,q),\quad \forall(p,q)\in P\odot P.
\end{equation}

\begin{Prop}\label{prop-prophi}
Given a right principal $\calG$-bundle, the division map $\phi_P$ from $P\odot P$ to $\calG$ has the following properties:
\begin{itemize}
\item[i)] for any point $(p,q)$ of $P\odot P$, we have
\[
\phi_P(p,q)\in \calG_{\varepsilon(q),\varepsilon(p)}.
\] 
\item[ii)] On the diagonal submanifold of the total space of $P\odot P$, we have
\[
\phi_P(p,p)=\iota_{\calG}(\varepsilon(p)),\quad \forall p\in P.
\] 
\item[iii)] for any pair $(p,q)\in P\odot P$, the following equation holds
\[
\phi_P(p,q)=\phi_P(q,p)^{-1};
\] 
notice that the previous equation makes sense, since $(p,q)\in P\odot P$ implies that $(q,p)\in P\odot P$ also.
\item[iv)] The triple $\left(\phi_P,\id_{\calG^2},\id_{X_{\calG}^2}\right)$ is an equivariant map from the right $\calG^2$-space $P\odot P$ to the right $\calG^2$-space $\calG$ endowed with the right generalized conjugation defined in Remark~\ref{rem-conjright} by the pair of maps $\left(J_{\conj}^R,\Psi_{\conj}^R\right)$ in Subsubsection~\ref{sssec-genconjgroupoid}.  
\end{itemize}
\end{Prop}
\begin{proof}
\begin{itemize}
\item[i)] As already seen, for any pair $(p,q)$ in $P\odot P$, one has
\[
t_{\calG}\left(\phi_P(p,q)\right)=\varepsilon(p).
\] 
On the other hand, Equation (\ref{eq-defphi}) implies, since $\calG$ acts from the right on $P$, that
\[
s_{\calG}\left(\phi_P(p,q)\right)=\varepsilon(q),
\]
whence the claim follows.
\item[ii)] Again, I make use of Equation (\ref{eq-defphi}): namely, for any pair $(p,p)$ it implies
\[
p=p\phi_P(p,p),
\] 
whence it follows, by Remark~\ref{rem-freeness},
\[
\phi_P(p,p)=\iota_{\calG}\left(\varepsilon(p)\right).
\]
\item[iii)] This follows immediately from Equation (\ref{eq-defphi}):
\[
q=p\phi_P(p,q)\Leftrightarrow p=q\phi_P(q,p),\quad \forall (p,q)\in P\odot P.
\] 
\item[iv)] First of all, one has to show the commutativity of the first diagram in Definition~\ref{def-equivgroupoid}; recall that, in this context, $M$ is the fibred product bundle $P\odot P$, $N$ is the manifold of arrows $\calG$, the smooth map $J_M$ is the product $\varepsilon\times\varepsilon$ and $J_N$ is $J_{\conj}$ from Proposition~\ref{prop-genconj}, and $\Theta$ is $\phi_P$ and $\varphi$ is the identity of $X_{\calG}^2$.

Thus, one can compute directly, with the help of the result proved in $i)$:
\begin{align*} 
\left(J_{\conj}\circ \phi_P\right)\!(p,q)&=J_{\conj}\!\left(\phi_P(p,q)\right)=\\
&=(\left(t_{\calG}\left(\phi_P(p,q)\right),s_{\calG}\left(\phi_P(p,q)\right)\right)=\\
&=\left(\varepsilon(p),\varepsilon(q)\right)=\\
&=\left(\id_{X_{\calG}^2}\circ (\varepsilon\times\varepsilon)\right)\!(p,q),\quad \forall (p,q)\in P\odot P.
\end{align*}

To prove the commutativity of the second diagram, consider a general $4$-tuple $(p,q;g_1,g_2)$ in $\left(P\odot P\right)\times_{\varepsilon\times\varepsilon}\left(\calG^2\right)$; then one has by Equation (\ref{eq-defphi})
\begin{align*}
qg_2&=pg_1 \phi_P(pg_1,qg_2)=\\
&=p\phi_P(p,q)g_2. 
\end{align*}
Notice that all identities make sense, in virtue of the commutativity of the commutativity of the first diagram in Definition~\ref{def-equivgroupoid}.

The freeness of the right $\calG$-action on $P$ from Remark~\ref{rem-freeness} implies that
\[
\phi_P(pg_1,q g_2)=g_1^{-1}\phi_P(p,q)g_2,
\]  
which is equivalent to the identity
\begin{align*}
&\left(\Psi_{\conj}\circ \left(\phi_P\times \id_{\calG^2}\right)\right)\!(p,q;g_1,g_2)=\left(\phi_P\circ\Psi_{P\odot P}\right)\!(p,q;g_1,g_2),\\&\forall (p,q;g_1,g_2)\in \left(P\odot P\right)\times_{\varepsilon\times\varepsilon}\left(\calG^2\right),
\end{align*}
where $\Psi_{P\odot P}$ denotes the right action map for the right $\calG^2$-space $P\odot P$.
\end{itemize}
\end{proof}
Let me end this subsection by giving the division map of the product bundle $P1\times P_2$, with $P_1$ and $P_2$ as in the hypotheses of Lemma~\ref{lem-prodbun} (from which one can easily deduce the division map for the fibred product bundle of $P$ and $\widetilde{P}$, two principal $\calG$-bundles over the same base space):
\begin{equation*}
\phi_{P_1\times P_2}\!\left((p_1,p_2),\left(\tildep_1,\tildep_2\right)\right)=\left(\phi_{P_1}(p_1,\tildep_1),\phi_{P_2}(p_2,\tildep_2)\right),
\end{equation*}
where $\pi_1(p_1)=\pi_1(\tildep_1)$ and $\pi_2(p_2)=\pi_2(\tildep_2)$.

\section{Equivariant morphisms between principal bundles and generalized gauge transformations}\label{sec-eqmorphgengauge}
In Section~\ref{sec-princbun}, I defined the notion of principal bundles with structure groupoid and discussed the notion of fibred product bundle of any two principal bundles; finally, I associated to any principal bundle $P$ with groupoid structure a canonical (twisted) equivariant map from the fibred product bundle of $P$ with itself to the structure groupoid itself.

In this Section, $i)$ I first review the concept of fibre-preserving, bundle morphisms between principal bundles with the same groupoid structure and over the same base manifold $M$, and $ii)$ I develop a theory of generalized gauge transformations in the sense of~\cite{C1}; the main tools for the development of such a theory are the notion of fibred product bundle and the canonical division map from Proposition~\ref{prop-prophi}.

\subsection{Equivariant maps between principal bundles}\label{ssec-eqmorph}  
I consider any two principal bundles $P_1$, $P_2$, with the same structure groupoid $\calG$ and over the same base manifold $M$.
\begin{Def}\label{def-equivbun}
A fibre-preserving, $\calG$-equivariant map between the principal bundles $P_1$ and $P_2$ (shortly, a bundle morphism between $P_1$ and $P_2$) is a twisted equivariant map $\left(\sigma,\id_{\calG},\id_{X_{\calG}}\right)$ in the sense of Definition~\ref{def-equivgroupoid} of Subsection~\ref{ssec-equivgroupoid} from the right $\calG$-space $P_1$ to the right $\calG$-space $P_2$, with the additional property to be {\em fibre-preserving} in the following sense:
\[
\pi_1\circ \sigma=\pi_2.
\]
(Notice that in particular $\sigma$ is also momentum-preserving.)
\end{Def}
It is not difficult to check that, for a triple $(P_1,P_2,P_3)$ of principal $\calG$-bundles over the base manifold $M$ and bundle morphisms $\sigma_{12}$ from $P_1$ to $P_2$ and $\sigma_{23}$ from $P_2$ to $P_3$, their composition is, by its very definition, again a bundle morphism from $P_1$ to $P_3$.
Clearly, the identity map $\id_P$ of $P$ is a bundle morphism from the principal $\calG$-bundle $P$ to itself.

Thus, it makes sense introduce the category $\mathsf{Bun}_{\calG,M}$ by the following assignments:
\begin{itemize}
\item[i)] {\bf Objects}: the objects of $\mathsf{Bun}_{\calG,M}$ are the principal $\calG$-bundles over the base manifold $M$
\item[ii)] {\bf Morphisms}: a morphism between two objects $P_1$, $P_2$ the category $\mathsf{Bun}_{\calG,M}$ is a bundle morphism from $P_1$ to $P_2$ in the sense of Definition~\ref{def-equivbun} 
\end{itemize}
Morphisms in the category $\mathsf{Bun}_{\calG,M}$ have the remarkable property of being bijective, as the following Lemma shows
\begin{Lem}\label{lem-inverequiv}
Every morphism $\sigma$ in $\mathsf{Mor}_{\mathsf{Bun}_{\calG,M}}\!(P_1,P_2)$, for any two objects $P_1$, $P_2$ of $\mathsf{Bun}_{\calG,M}$, is bijective.
\end{Lem}
\begin{proof}
One has to show $i)$ injectivity and $ii)$ surjectivity of $\sigma$.
Let me first show injectivity.
Namely, consider two points $p_1$, $q_1$ of $P_1$, such that
\[
\sigma(p_1)=\sigma(q_1).
\]
Since $\sigma$ is fibre-preserving, it follows that $p_1$ and $q_1$ lie in the same fibre, whence
\[
q_1=p_1 \phi_{P_1}(p_1,q_1).
\] 
The $\calG$-equivariance of $\sigma$ implies readily
\[
\sigma(q_1)=\sigma\left(p_1\phi_{P_1}(p_1,q_1)\right)=\sigma(p_1)\phi_{P_1}(p_1,q_1)\overset{!}=\sigma(p_1).
\]
The freeness of the action of $\calG$ on $P_2$ implies
\[
\phi_{P_1}(p_1,q_1)=\iota_{\calG}\!\left(\varepsilon_2(\sigma(p_1))\right)=\iota_{\calG}\left(\varepsilon_1(p_1)\right),
\]
whence
\[
q_1=p_1\iota_{\calG}\left(\varepsilon_1(p_1)\right)=p_1.
\]

As for surjectivity, one has to show that for any point $p_2$ of $P_2$, there exists a point $p_1$ of $P_1$, such that
\[
\sigma(p_1)=p_2.
\]
From the fact that $\sigma$ is fibre-preserving, it follows immediately that surjectivity is a fibrewise property for equivariant morphisms.
Consider therefore a point $p_2$ of $P_2$ and we take its projection $\pi_2(p_2)=\colon x$; by the surjectivity of $\pi_1$, consider a point $q_{1,x}$ of $P_1$, such that $\pi_1(q_{1,x})=x$.
Consider further the image w.r.t.\ $\sigma$ of $q_{1,x}$; it lies in the same fibre of $p_2$, whence
\[
p_2=\sigma(q_{1,x})\phi_{P_2}(q_{1,x},p_2).
\]
By Proposition~\ref{prop-prophi} and Definition~\ref{def-equivbun},
\[
t_{\calG}\!\left(\phi_{P_2}(q_{1,x},p_2)\right)=\varepsilon_2\!\left(\sigma(q_{1,x})\right)=\varepsilon_1(q_{1,x}),
\]
hence, one can form the element
\[
p_1\colon=q_{1,x}\phi_{P_2}(\sigma(q_{1,x}),p_2).
\]
An easy computation, using the $\calG$-equivariance of $\sigma$, gives
\begin{align*}
\sigma(p_1)&=\sigma\left(q_{1,x}\phi_{P_2}(\sigma(q_{1,x}),p_2)\right)=\\
&=\sigma(q_{1,x})\phi_{P_2}(\sigma(q_{1,x}),p_2)=\\
&=p_2,
\end{align*}
by Equation (\ref{eq-defphi}).
\end{proof}

\begin{Rem}
M{\oe}rdijk showed that every bundle morphism between principal $\calG$-bundles is an isomorphism, by reducing the problem to {\em trivial principal bundles}.
Later, I will give another characterisation of morphisms of the category $\mathsf{Bun}_{\calG,M}$ and I will also see that, in fact, every morphism is an isomorphism, by using global arguments.
\end{Rem}

\subsection{Generalized gauge transformations}\label{ssec-gengauge}
I want now to discuss a different characterization of bundle morphisms between principal bundles with structure groupoid; in the previous Subsection, we have viewed bundle morphisms between principal $\calG$-bundles as special types of equivariant morphisms between right $\calG$-spaces in the sense of Definition~\ref{def-equivgroupoid}.

I define now generalized gauge transformations between two principal $\calG$-bundles analogously to what I did in Section 4 of~\cite{C1}, although the fact that I deal with groupoids, instead of groups, requires more care; but the idea is nonetheless the same, i.e.\ to consider maps from the fibred product of two principal $\calG$-bundles to the structure groupoid $\calG$ itself, satisfying some particular properties.
\begin{Def}\label{def-gengauge}
A {\em generalized gauge transformation between the principal $\calG$-bundles $P_1$ and $P_2$} is, by definition, a (twisted) equivariant map $\left(K,\id_{\calG^2},\id_{X_{\calG}^2}\right)$ from the right $\calG^2$-space $P_1\odot P_2$ and $\calG$, viewed as a right $\calG^2$-space via the right generalized conjugation defined in Remark~\ref{rem-conjright} via the pair of maps $\left(\overline{J}_{\conj}^R,\overline{\Psi}_{\conj}^R\right)$ in Subsubsection~\ref{sssec-genconjgroupoid}.
\end{Def}
The set of all generalized gauge transformations between the principal $\calG$-bundles $P_1$ and $P_2$ is denoted by $C^{\infty}\!\left(P_1\odot P_2,\calG\right)^{\calG^2}$.

\begin{Rem}\label{rem-explgengauge}
Let me give a more detailed account of the properties of generalized gauge transformations.
First of all, a generalized gauge transformation between the principal $\calG$-bundles $P_1$ and $P_2$ is a smooth map from the fibred product $P_1\odot P_2$ to the structure groupoid $\calG$.
The fact that the triple $\left(K,\id_{\calG^2},\id_{X_{\calG}^2}\right)$ is a twisted bundle morphism from $P_1\odot P_2$ to $\calG$, viewed both as right $\calG^2$-spaces, can be translated into the following set of equations:
\begin{itemize}
\item[i)] The first diagram of Definition~\ref{def-equivgroupoid} implies immediately
\[
\overline{J}_{\conj}^R\circ K=\varepsilon_1\times\varepsilon_2\Rightarrow \begin{cases}
s_{\calG}\!\left(K(p_1,p_2)\right)&=\varepsilon_1(p_1),\\
t_{\calG}\!\left(K(p_1,p_2)\right)&=\varepsilon_2(p_2),\quad \forall(p_1,p_2)\in P_1\odot P_2.
\end{cases}
\] 
\item[ii)] The second diagram of Definition~\ref{def-equivgroupoid}, which defines precisely the $\calG^2$-equivariance, may be restated as follows:
\[
K(p_1 g_1,p_2 g_2)=g_2^{-1} K(p_1,p_2) g_1,\quad (p_1,p_2;g_1,g_2)\in \left(P_1\odot P_2\right)\times_{\varepsilon_1\times \varepsilon_2}\left(\calG^2\right).
\] 
\end{itemize}
\end{Rem}

I now prove the following
\begin{Thm}\label{thm-gengaugeeq}
The set $\mathsf{Mor}_{\mathsf{Bun}_{\calG,M}}\!(P_1,P_2)$ of morphisms between principal $\calG$-bundles $P_1$ and $P_2$ in the sense of Definition~\ref{def-equivbun} is in one-to-one correspondence to the set $C^{\infty}\!\left(P_1\odot P_2,\calG\right)^{\calG^2}$ of generalized gauge transformations between $P_1$ and $P_2$.
\end{Thm}
\begin{proof}
Consider first a bundle morphism $\sigma$ between $P_1$ and $P_2$ in the sense of Definition~\ref{def-equivbun}; I associate to it the following composite map
\[
\sigma\leadsto K_\sigma\!(p_1,p_2)\colon=\phi_{P_2}\!\left(p_2,\sigma(p_1)\right),\quad \forall (p_1,p_2)\in P_1\odot P_2.
\]
First of all, $K_\sigma$ is a well-defined map from $P_1\odot P_2$ to $\calG$: namely, it is already known from Section~\ref{sec-princbun} that the map $\phi_{P_2}$ is a smooth map from $P_2\odot P_2$ to $\calG$, and, if the pair $(p_1,p_2)$ belongs to the fibred product $P_1\odot P_2$, then the pair $(p_2,\sigma(p_1))$ belongs to the fibred product of $P_2$ with itself, since
\[
\pi_2(\sigma(p_1))=\pi_1(p_1)=\pi_2(p_2).
\]
It remains to show the commutativity of the diagrams in Definition~\ref{def-equivgroupoid}, which have been translated in two sets of equations in Remark~\ref{rem-explgengauge}.
Using the properties of the map $\phi_{P_2}$, displayed in Proposition~\ref{prop-prophi}, one shows commutativity of the first diagram, namely:
\begin{align*}
s_{\calG}\!\left(K_\sigma(p_1,p_2)\right)&=s_{\calG}\!\left(\phi_{P_2}(p_2,\sigma(p_1))\right)=\\
&=\varepsilon_2(\sigma(p_1))=\\
&=\varepsilon_1(p_1),
\end{align*}
and
\begin{align*}
t_{\calG}\!\left(K_\sigma(p_1,p_2)\right)&=t_{\calG}\!\left(\phi_{P_2}(p_2,\sigma(p_1))\right)=\\
&=\varepsilon_2(p_2).
\end{align*}

As for the second diagram, one gets, again by Proposition~\ref{prop-prophi},
\begin{align*}
K_\sigma(p_1 g_1,p_2 g_2)&=\phi_{P_2}\!\left(p_2 g_2,\sigma(p_1 g_1)\right)=\\
&\overset{\text{by $\calG$-equivariance of $\sigma$}}=\phi_{P_2}\!\left(p_2 g_2,\sigma(p_1) g_1\right)=\\
&=g_2^{-1} \phi_{P_2}\!\left(p_2,\sigma(p_1)\right) g_1=\\
&= g_2^{-1} K_{\sigma}(p_1,p_2) g_1,\quad \forall (p_1,p_2;g_1,g_2)\in \left(P_1\odot P_2\right)\times_{\varepsilon_1\times\varepsilon_2}\left(\calG^2\right).
\end{align*}

On the other hand, given a generalized gauge transformation $K$ between $P_1$ and $P_2$, it is possible to define a bundle morphism $\sigma_K$ from $P_1$ to $P_2$ by the following rule:
\[
\sigma_K(p_1)\colon=p_2 K(p_1,p_2),\quad (p_1,p_2)\in P_1\odot P_2.
\]
The previous formula is well-defined, in the following sense: $i)$ the right multiplication makes sense and $ii)$ it does {\em not} depend on the choice of $p_2$, as long as the pair $(p_1,p_2)$ belongs to $P_1\odot P_2$.
To prove $i)$, notice that
\[
t_{\calG}\!\left(K(p_1,p_2)\right)=\varepsilon_2(p_2)\Rightarrow \left(p_2,K(p_1,p_2)\right)\in P_2\times_{\varepsilon_2}\calG,\forall p_1\in P_1\quad\text{s.t.}\quad (p_1,p_2)\in P_1\odot P_2.  
\]
To prove $ii)$, consider, for $p_1$ in $P_1$ fixed, another pair $(p_1,q_2)$ in $P_1\odot P_2$; it follows immediately, by Definition~\ref{def-princgroupoid}, that 
\[
q_2=p_2 \phi_{P_2}(p_2,q_2),
\]
whence
\begin{align*}
\sigma_K(p_1)&=q_2 K(p_1,q_2)=\\
&=p_2 \phi_{P_2}(p_2,q_2) K\!\left(p_1,p_2\phi_{P_2}(p_2,q_2)\right)=\\
&=p_2 \phi_{P_2}(p_2,q_2)\phi_{P_2}(p_2,q_2)^{-1}K(p_1,p_2)=\\
&=p_2 t_{\calG}\!\left(\phi_{P_2}(p_2,q_2)\right)K(p_1,p_2)=\\
&=p_2 \varepsilon_2(p_2) K(p_1,p_2)=\\
&=p_2 K(p_1,p_2).
\end{align*}

It remains to show that the triple $\left(\sigma_K,\id_{\calG^2},\id_{X_{\calG}^2}\right)$ is a bundle morphism between $P_1$ and $P_2$; this is equivalent to showing the commutativity of the two diagrams in Definition~\ref{def-equivgroupoid}.
To show the commutativity of the first diagram, we compute
\begin{align*}
\varepsilon_2\!\left(\sigma_K(p_1)\right)&=\varepsilon_2\!\left(p_2 K(p_1,p_2)\right)=\\
&=s_{\calG}\!\left(K(p_1,p_2)\right)=\\
&=\varepsilon_1(p_1),\quad \forall p_1\in P_1, p_2\quad\text{s.t.}\quad (p_1,p_2)\in P_1\odot P_2.
\end{align*}
The commutativity of the second diagram follows by the following computation:
\begin{align*}
\sigma_K(p_1 g_1)&=p_2 K(p_1 g_1,p_2)=\\
&=p_2 K(p_1,p_2) g_1=\\
&\sigma_K(p_1) g_1,\quad (p_1,g_1)\in P_1\times_{\varepsilon_1}\calG.
\end{align*}
The property of $\sigma_K$ being fibre-preserving follows from
\begin{align*}
\pi_2\!\left(\sigma_K(p_1)\right)&=\pi_2\!\left(p_2 K(p_1,p_2)\right)=\\
&=\pi_2(p_2)=\\
&=\pi_1(p_1),\quad p_1\in P_1,
\end{align*}
since the pair $(p_1,p_2)$ belongs to $P_1\odot P_2$.

One has to show that the assignments
\[
\sigma\leadsto K_\sigma\quad\text{and}\quad K\leadsto \sigma_K
\]
are inverse to each other.
A direct computation shows
\begin{align*}
\sigma_{K_{\sigma}}(p_1)&=p_2 K_{\sigma}(p_1,p_2)=\\
&=p_2 \phi_{P_2}\!\left(p_2,\sigma(p_1)\right)=\\
&=\sigma(p_1),\quad \forall p_1\in P_1;
\end{align*}
on the other hand,
\begin{align*}
K_{\sigma_K}(p_1,p_2)&=\phi_{P_2}\!\left(p_2,\sigma_K(p_1)\right)=\\
&=\phi_{P_2}\!\left(p_2,p_2 K(p_1,p_2)\right)=\\
&=\phi_{P_2}(p_2,p_2) K(p_1,p_2)=\\
&=\iota_{\calG}(\varepsilon_2(p_2))K(p_1,p_2)=\\
&=K(p_1,p_2),\quad \forall (p_1,p_2)\in P_1\odot P_2.
\end{align*}
\end{proof}

It was proved in Lemma~\ref{lem-inverequiv} of Subsection~\ref{ssec-eqmorph} that any bundle morphism between any two $\calG$-principal bundles $P_1$ and $P_2$ is bijective.
Moreover, every bundle morphism $\sigma$ between $P_1$ and $P_2$ is invertible: namely, consider the generalized gauge transformation $K_\sigma$, canonically associated to $\sigma$ by Theorem~\ref{thm-gengaugeeq}.

\begin{Lem}\label{lem-gengaugeinv1}
For any bundle morphism $\sigma$ between $P_1$ and $P_2$, the map
\[
P_2\odot P_1\ni(p_2,p_1)\overset{K_{\sigma^{-1}}}\mapsto K_\sigma(p_1,p_2)^{-1}\in \calG, 
\]
defines a generalized gauge transformation between $P_2$ and $P_1$.
\end{Lem}
\begin{proof}
First of all, notice that the definition makes sense, since 
\[
(p_1,p_2)\in P_1\odot P_2\Leftrightarrow(p_2,p_1)\in P_2\odot P_1. 
\]
It remains to show the commutativity of the two diagrams in Definition~\ref{def-equivgroupoid}.
To show the commutativity of the first one, one computes
\begin{align*}
\overline{J}_{\conj}^R\!\left(K_{\sigma^{-1}}\!(p_2,p_1)\right)&=\overline{J}_{\conj}^R\!\left(K_\sigma(p_1,p_2)^{-1}\right)=\\
&=\left(s_{\calG}\!\left(K_\sigma(p_1,p_2)^{-1}\right),t_{\calG}\!\left(K_\sigma(p_1,p_2)^{-1}\right)\right)=\\
&=\left(t_{\calG}\!\left(K_\sigma(p_1,p_2)\right),s_{\calG}\!\left(K_\sigma(p_1,p_2)\right)\right)=\\
&=\left(\varepsilon_2(p_2),\varepsilon_1(p_1)\right)=\\
&=\left(\varepsilon_2\times\varepsilon_1\right)\!(p_2,p_1),\quad \forall (p_2,p_1)\in P_2\odot P_1.
\end{align*}
The commutativity of the second diagram follows from 
\begin{align*}
K_{\sigma^{-1}}\!(p_2 g_2,p_1 g_1)&=K_\sigma(p_1 g_1,p_2 g_2)^{-1}=\\
&=\left(g_2^{-1}K_\sigma(p_1,p_2) g_1\right)=\\
&=g_1^{-1} K_\sigma(p_1,p_2) g_2=\\
&=g_1^{-1} K_{\sigma^{-1}}\!(p_2,p_1) g_2,\quad \forall (p_2,p_1;g_2,g_1)\in \left(P_2\odot P_1\right)\times_{\varepsilon_2\times \varepsilon_1}\left(\calG^2\right).
\end{align*}
\end{proof}
Hence, to any bundle morphism $\sigma$ between $P_1$ and $P_2$ one associates in a canonical way two generalized gauge transformations, $K_\sigma$ between $P_1$ and $P_2$, and $K_{\sigma^{-1}}$ between $P_2$ and $P_1$.

The next lemma shows their relationship explicitly.
\begin{Lem}\label{lem-gengaugeinv2}
The unique bundle morphism $\tau$ associated to the generalized gauge transformation $K_{\sigma^{-1}}$, for any bundle morphism $\sigma$ from $P_1$ to $P_2$, by Theorem~\ref{thm-gengaugeeq}, is the inverse map to $\sigma$.
\end{Lem}
\begin{proof}
By definition the bundle morphism $\tau$ satisfies the equation
\[
\tau(p_2)=p_1 K_{\sigma^{-1}}\!(p_2,p_1),\quad \forall p_2\in P_2,
\]
and the pair $(p_2,p_1)$ belongs to $P_2\odot P_1$.
Then, by a direct computation on gets:
\begin{align*}
\sigma(\tau(p_2))&=\sigma\left(p_1K_{\sigma^{-1}}\!(p_2,p_1)\right)=\\
&=\sigma(p_1)K_{\sigma^{-1}}\!(p_2,p_1)=\\
&=\sigma(p_1) K_\sigma(p_1,p_2)^{-1}=\\
&=\sigma(p_1) \phi_{P_2}\!\left(p_2,\sigma(p_1)\right)^{-1}=\\
&=\sigma(p_1)\phi_{P_2}\!\left(\sigma(p_1),p_2\right)=\\
&=p_2,\quad \forall p_2\in P_2,
\end{align*}
where the pair $(p_2,p_1)$ belongs to $P_2\odot P_1$.

On the other hand, one has
\begin{align*}
\tau(\sigma(p_1))&=p_1 K_{\sigma^{-1}}\!(\sigma(p_1),p_1)=\\
&=p_1 K_\sigma(p_1,\sigma(p_1))^{-1}=\\
&=p_1 \phi_{P_2}\!\left(\sigma(p_1),\sigma(p_1)\right)^{-1}=\\
&=p_1 \phi_{P_2}(\sigma(p_1),\sigma(p_1))=\\
&=p_1 \iota_{\calG}(\varepsilon_2(\sigma(p_1)))=\\
&=p_1 \iota_{\calG}(\varepsilon_1(p_1))=\\
&=p_1,
\end{align*}
where was used the fact that the pair $(p_1,\sigma(p_1))$ belongs to $P_1\odot P_2$.
\end{proof}

\subsubsection{Gauge transformations}\label{sssec-gaugetrsf}
In this Subsubsection, I want to study the notion of gauge transformation of a principal $\calG$-bundle $P$.
First of all, I consider a bundle morphism $\sigma$ from $P$ to itself in the sense of Definition~\ref{def-equivbun}.
By Theorem~\ref{thm-gengaugeeq}, there is a unique generalized gauge transformation $K_\sigma$ on $P$, i.e.\ a $\calG^2$-equivariant map from the fibred product of $P$ with itself, defined via
\[
K_\sigma(p,q)\colon=\phi_P(q,\sigma(p)),
\]
where $\phi_P$ is the canonical map associated to the bundle $P$, thoroughly discussed in Subsection~\ref{ssec-propphi}.
On the other hand, since $\sigma$ is fibre-preserving, one has
\[
\sigma(p)=p G_\sigma(p),\quad\forall p\in P,
\]
for a unique element $G_\sigma(p)$, depending smoothly on $P$ and belonging to the groupoid $\calG$.
By the freeness of the action of $\calG$ on $P$, it follows
\[
G_\sigma(p)=\phi_P(p,\sigma(p))=K_\sigma(p,p),
\]
i.e.\ $G_\sigma$ is the restriction to the diagonal of $P\odot P$ of the unique generalized gauge transformation associated to $\sigma$.

Moreover, by the properties of generalized gauge transformations, 
\[
\begin{cases}
s_{\calG}\!\left(G_\sigma(p)\right)&=s_{\calG}\!\left(K_\sigma(p,p)\right)=\varepsilon(p),\\
t_{\calG}\!\left(G_\sigma(p)\right)&=t_{\calG}\!\left(K_\sigma(p,p)\right)=\varepsilon(p),
\end{cases}
\]
i.e.\ $g_\sigma(p)\in \calG_{\varepsilon(p),\varepsilon(p)}$, for any $p\in P$.
Furthermore, the following equivariance property of $G_\sigma$ holds:
\begin{align*}
G_\sigma(p g)&=\phi_P\!\left(pg,\sigma(pg)\right)=\\
&=\phi_P\!\left(pg,\sigma(p)g\right)=\\
&=g^{-1} \phi_P\!\left(p,\sigma(p)\right) g=\\
&=g^{-1} G_\sigma(p) g,\quad \forall (p,g)\in P\times_\varepsilon \calG.
\end{align*}

On the other hand, if we have a smooth map from $P$ to $\calG$, satisfying
\begin{align}
\label{eq-group}&G(p)\in\calG_{\varepsilon(p),\varepsilon(p)},\quad\forall p\in P,\\
\label{eq-equivgauge}&G(p g)=g^{-1} G(p) g,\quad (p,g)\in P\times_\varepsilon \calG,
\end{align}
the well-defined assignment 
\[
p\overset{\sigma_K}\mapsto p G(p),\quad \forall p\in P,
\]
defines in an obvious way a bundle morphism on $P$, which I denote by $\sigma_G$.
Therefore, by Theorem~\ref{thm-gengaugeeq}, $\sigma_G$ defines a unique generalized gauge transformation $K_G$ by the rule
\begin{equation}\label{eq-gaugegengauge}
\begin{aligned}
K_G(p,q)&\colon=\phi_P\!\left(q,\sigma_G(p)\right)=\\
&=\phi_P\!\left(q,p G(p)\right)=\\
&=\phi_P\!\left(q,p\right)G(p)=\\
&=\phi_P(p,q)^{-1}G(p),\quad \forall (p,q)\in P\odot P.
\end{aligned}
\end{equation}
Computations similar to those used in the proof of Theorem~\ref{thm-gengaugeeq} imply that the assignments
\[
G\leadsto K_G\quad\text\quad K\leadsto\iota_\Delta^*\!K,
\]
are inverse to each other, where $G$ is any map from $P$ to $\calG$, satisfying both Equations (\ref{eq-group}) and (\ref{eq-equivgauge}), and $K$ is a generalized gauge transformation of $P$; the map $\iota_\Delta$ denotes here the imbedding of the diagonal of $P\odot P$.

By Equation (\ref{eq-group}), one can define on the set $C^{\infty}(P,\calG)^\calG$ of maps from $P$ to $\calG$ satisfying Equations (\ref{eq-group}) and (\ref{eq-equivgauge}), a product structure: in fact,
\[
(G_1G_2)(p)\colon=G_1(p)G_2(p),\quad \forall p\in P. 
\]
It is clear that the map $G_1 G_2$ enjoys again Equations (\ref{eq-group}) and (\ref{eq-equivgauge}).
Moreover, the product is associative, since, for any $p\in P$, the set $\calG_{\varepsilon(p),\varepsilon(p)}$ is a group.

Furthermore, the map
\[
\iota_{\calG}\circ\varepsilon\colon P\to \calG  
\]
satisfies both Equations (\ref{eq-group}) and (\ref{eq-equivgauge}): in fact, e.g.\ Equation (\ref{eq-equivgauge}) holds because:
\begin{align*}
\iota_{\calG}\!\left(\varepsilon(pg)\right)&=\iota_{\calG}\!\left(s_{\calG}(g)\right)=\\
&=g^{-1}g=\\
&=g^{-1} \iota_{\calG}\!\left(\varepsilon(p)\right)g,\quad \forall (p,g)\in P\times_\varepsilon \calG.
\end{align*}
It is not difficult to prove that $\iota_{\calG}\circ\varepsilon$ is the unit for the product in $C^{\infty}(P,\calG)^\calG$.
Analogously, to any $G$ in $C^{\infty}(P,\calG)^\calG$, the map
\[
G^{-1}\!(p)\colon=G(p)^{-1},\quad \forall p\in P,
\]
belongs again to $C^{\infty}(P,\calG)^\calG$, and a direct computation shows that
\[
G G^{-1}=G^{-1} G=\iota_{\calG}\circ\varepsilon.
\]
Hence, it is possible to summarize all the computations so far in the following
\begin{Prop}\label{prop-gaugegr}
For any principal $\calG$-bundle $P$, the set $C^{\infty}(P,\calG)^\calG$ of maps from $P$ to $\calG$, satisfying Equations (\ref{eq-group}) and (\ref{eq-equivgauge}), is in one-to-one correspondence via the maps
\[
C^{\infty}(P,\calG)^\calG\ni G\leadsto K_G\in C^{\infty}(P\odot P,\calG)^{\calG^2}\ni K\leadsto \iota_\Delta^*\!K\in C^{\infty}(P,\calG)^\calG
\]
with the set of generalized gauge transformations $C^{\infty}(P\odot P,\calG)^{\calG^2}$; here, $K_G$ denotes the map defined by Equation (\ref{eq-gaugegengauge}).
Moreover, the set $C^\infty(P,\calG)^\calG$ is a group, called the {\em gauge group of $P$}; thus, the set of bundle (auto)morphisms of $P$, being in one-to-one correspondence with the gauge group $C^\infty(P,\calG)^\calG$, inherits a group structure via composition, and the map $\sigma\mapsto G_\sigma$, for any bundle (auto)morphism $\sigma$ of $P$, is an isomorphism of groups.
\end{Prop}
\begin{Rem}
The previous Proposition implies readily that the gauge group of a principal $\calG$-bundle $P$ may be viewed as the isotropy group at $P$ of the groupoid of generalized gauge transformations, which I will introduce and discuss in the Subsection~\ref{ssec-gengaugegroupoid}.
\end{Rem}

\subsection{Invariance of the division map w.r.t.\ bundle isomorphisms}\label{ssec-indvivmap}
In this short subsection, I will display a trivial, but important property of the division map, namely its invariance w.r.t.\ bundle morphisms.
In other words: one already knows that bundle morphisms between right principal $\calG$-bundles over the same base space are isomorphisms.
Therefore, considering isomorphism classes of $\calG$-bundles over the same base space, one may consider one representative $P$, and consider subsequently its division map $\phi_P$: the latter is an invariant of the isomorphism class, i.e.\ it does not depend on the choice of the representative.

More formally, the content of the previous discussion may be restated in the following
\begin{Thm}\label{thm-gaugeinvdiv}
Let $\sigma$ be a bundle morphism from the right principal $\calG$-bundle $P_1$ to the right principal $\calG$-bundle $P_2$, both over the same base space.

Then, the following identity holds:
\begin{equation}\label{eq-gaugeinvdiv}
\phi_{P_2}\circ\!\left(\sigma\times \sigma\right)=\phi_{P_1}\quad\text{on $P_1\odot P_1$}.
\end{equation}
\end{Thm} 
\begin{proof}
First of all, let us check that the map on the left-hand side is well-defined.
This is not difficult: in fact, considering a pair $(p_1,\overline{p}_1)$ in $P_1\odot P_1$, it follows immediately that the pair $\left(\sigma(p_1),\sigma(\overline{p}_1)\right)$ belongs to $P_2\odot P_2$, since $\sigma$ is fibre-preserving.

Second, the identity follows from the following computation:
\begin{align*}
\sigma(p_1)\phi_{P_2}\!\left(\sigma(p_1),\sigma(\overline{p}_1)\right)&=\sigma(\overline{p}_1)=\\
&=\sigma\!\left(p_1\phi_{P_1}(p_1,\overline{p}_1)\right)=\\
&=\sigma(p_1)\phi_{P_1}(p_1,\overline{p}_1),\quad \forall (p_1,\overline{p}_1)\in P_1\odot P_1,
\end{align*}
by the definitions of the division maps $\phi_{P_1}$ and $\phi_{P_2}$, and by the $\calG$-equivariance of $\sigma$.
Since the action of $\calG$ is free, Identity (\ref{eq-gaugeinvdiv}) follows immediately.  
\end{proof}

As a simple consequence, the assignment to a bundle morphism $\sigma$ between $P_1$ and $P_2$, right principal $\calG$-bundles over the same base space, of a generalized gauge transformation of Theorem~\ref{thm-gengaugeeq} may be also rewritten as follows:
\[
\sigma\leadsto K_\sigma(p_1,p_2)=\phi_{P_1}\!\left(\sigma^{-1}(p_2),p_1\right),\quad (p_1,p_2)\in P_1\odot P_2.
\]

\subsection{The groupoid of generalized gauge transformations}\label{ssec-gengaugegroupoid}
As shown in Theorem~\ref{thm-gengaugeeq} in Subsection~\ref{ssec-gengauge}, any bundle morphism between two principal $\calG$-bundles is invertible, thus, by definition, every morphism of the category $\mathsf{Bun}_{\calG,M}$ is invertible, making $\mathsf{Bun}_{\calG,M}$ to an abstract groupoid. 
I want to discuss this groupoid from the point of view of generalized gauge transformations.
Let me begin with a notational remark:

\fbox{\parbox{12cm}{\bf A bundle morphism from an object $P_i$ to another object $P_j$ of the category $\mathsf{Bun}_{\calG,M}$ will be denoted by $\sigma_{ij}$}}

I consider now a triple $\left(P_1,P_2,P_3\right)$ of objects of $\mathsf{Bun}_{\calG,M}$, and corresponding bundle morphisms $\sigma_{12}$ and $\sigma_{23}$.
Since $\sigma_{23}\circ\sigma_{12}$ is obviously $G$-equivariant and fibre-preserving from $P_1$ to $P_3$ in the sense of Definition~\ref{def-equivbun}, there is a unique generalized gauge transformation in $C^\infty\!(P_1\odot P_3,\calG)^{\calG^2}$ by Theorem~\ref{thm-gengaugeeq}, given explicitly by:
\[
\left(\sigma_{23}\circ\sigma_{12}\right)\leadsto K_{\sigma_{23}\circ\sigma_{12}}\!\left(p_1,p_3\right)=\phi_{P_3}\!\left(p_3,\left(\sigma_{23}\circ\sigma_{12}\right)(p_1)\right),\quad \pi_1(p_1)=\pi_3(p_3).
\]
A direct computation shows
\begin{align*}
\left(\sigma_{23}\circ\sigma_{12}\right)(p_1)&=\sigma_{23}\left(p_2K_{\sigma_{12}}\!\left(p_1,p_2\right)\right)=\\
&=\sigma_{23}(p_2)K_{\sigma_{12}}\!\left(p_1,p_2\right)=\\
&=p_3 K_{\sigma_{23}}\!\left(p_2,p_3\right)K_{\sigma_{12}}\!\left(p_1,p_2\right),
\end{align*}
where $p_2\in P_2$ such that $\pi_2(p_2)=\pi_1(p_1)=\pi_3(p_3)$.
The freeness of the action of $\calG$ on $P_3$ implies finally
\begin{equation*}%\label{eq-composisom}
K_{13}\!\left(p_1,p_3\right)\colon=K_{\sigma_{23}\circ\sigma_{12}}\!\left(p_1,p_3\right)=K_{\sigma_{23}}\!\left(p_2,p_3\right)K_{\sigma_{12}}\!\left(p_1,p_2\right).
\end{equation*}
(In order to avoid cumbersome notations, I simply abbreviate $K_{\sigma_{12}}$ by $K_{12}$ and so on.)

%For any two objects $P_1$, $P_2$ of the category $\mathsf{Bun}_{\calG,M}$, we consider the set $C^{\infty}\!(P_1\odot P_2,\calG)^{\calG\times \calG}$ of {\em generalized gauge transformations from $P_1$ to $P_4$}.
%(We want to point out now that, in principle, the set $C^{\infty}\!(P_1\odot P_2,\calG)^{\calG\times \calG}$ may be empty, since $P_1$ and $P_2$ need not necessarily be isomorphic. 
%We shall hence assume formally any two objects $P_1$ and $P_2$ of $\mathsf{Bun}_{\calG,M}$ to be isomorphic; this assumption is not hollow, as we will provide later some examples.)

Now, for any triple $(P_1,P_2,P_3)$ of objects of $\mathsf{Bun}_{\calG,M}$, consider the product operation 
\begin{align*}
C^{\infty}\!(P_2\odot P_3,\calG)^{\calG^2}\times C^{\infty}\!(P_1\odot P_2,\calG)^{\calG^2}\ni\left(K_{23},K_{12}\right)\mapsto&\left(K_{23}\star K_{12}\right)(p_1,p_3)\colon=\\
&=K_{23}(p_2,p_3) K_{12}(p_1,p_2),
\end{align*}
for any pair $(p_1,p_3)$ in $P_1\odot P_3$, and $p_2$ in $P_2$ satisfying
\[
(p_1,p_2)\in P_1\odot P_2\Rightarrow (p_2,p_3)\in P_2\odot P_3.
\] 

First of all, the operation $\star$ is well-defined, since 
\[
s_{\calG}\!\left(K_{23}(p_2,p_3)\right)=\varepsilon_2(p_2)=t_{\calG}\!\left(K_{12}(p_1,p_2)\right).
\]
Moreover, since $K_{12}$ and $K_{23}$ are both generalized gauge transformations, their product $K_{23}\star K_{12}$ does not depend on the choice of $p_2\in P_2$, as long as $\pi_1(p_1)=\pi_2(p_2)=\pi_3(p_3)$ holds: namely, for another representative $q_2=p_2 \phi_{P_2}\!\left(p_2,q_2\right)$ in the same fibre of $\pi_2$, we get
\begin{align*}
K_{13}\!(p_1,p_3)&=K_{23}\!(p_2\phi_{P_2}\!\left(p_2,q_2\right),p_3)K_{12}\!(p_1,p_2 \phi_{P_2}\!\left(p_2,q_2\right))=\\
&=K_{23}\!(p_2,p_3)\phi_{P_2}\!\left(p_2,q_2\right) \phi_{P_2}\!\left(p_2,q_2\right)^{-1} K_{12}\!(p_1,p_2)=\\
&=K_{23}\!(p_2,p_3)\iota_{\calG}\!\left(t_{\calG}\!\left(\phi_{P_2}\!\left(p_2,q_2\right)\right)\right) K_{12}\!(p_1,p_2)=\\
&=K_{23}\!(p_2,p_3)\iota_{\calG}\!\left(\varepsilon_2(p_2)\right)K_{12}\!(p_1,p_2)=\\
&=K_{23}\!(p_2,p_3)K_{12}\!(p_1,p_2),
\end{align*} 
since
\[
s_{\calG}\!\left(K_{23}\!(p_2,p_3)\right)=\varepsilon_2(p_2).
\]
On the other hand, $K_{23}\star K_{12}$ is $\calG\times\calG$-equivariant:
\begin{align*}
K_{13}\!(p_1 g,p_3 h)&=K_{23}\!(p_2,p_3 h)K_{12}\!(p_1 g,p_2)=\\
&=h^{-1}K_{23}\!(p_2,p_3)K_{12}\!(p_1,p_2)g=\\
&=h^{-1}K_{13}\!(p_1 ,p_3) g,\quad \forall g,h\in \calG.
\end{align*}

Consider now a $4$-tuple $\left(P_1,P_2,P_3,P_4\right)$ of objects of the category $\mathsf{Bun}_{\calG,M}$, and the three respective sets of generalized gauge transformations:
\[
C^{\infty}\!(P_1\odot P_2,\calG)^{\calG^2},\quad C^{\infty}\!(P_2\odot P_3,\calG)^{\calG^2}\quad \text{and}\quad C^{\infty}\!(P_3\odot P_4,\calG)^{\calG^2}.
\]
It makes sense to consider the following iterated operations of the map $\star$:
\[
K_{34}\star\left(K_{23}\star K_{12}\right)\quad\text{and}\quad \left(K_{34}\star K_{23}\right)\star K_{12},
\]
for any $K_{12}\in C^{\infty}\!(P_1\odot P_2,\calG)^{\calG^2}$, $K_{23}\in C^{\infty}\!(P_2\odot P_3,\calG)^{\calG^2}$ and $K_{34}\in C^{\infty}\!(P_3\odot P_4,\calG)^{\calG^2}$.
Explicit computations give
\begin{align*}
\left(K_{34}\star\left(K_{23}\star K_{12}\right)\right)\!(p_1,p_4)&=K_{34}\!(p_3,p_4)\left(K_{23}\star K_{12}\right)\!(p_1,p_3)=\quad \text{$(\pi_3(p_3)=\pi_1(p_1))$}\\
&=K_{34}\!(p_3,p_4)K_{23}\!(p_2,p_3) K_{12}\!(p_1,p_2)=\quad \text{$(\pi_2(p_2)=\pi_1(p_1))$}\\
&=\left(K_{34}\star K_{23}\right)\!(p_2,p_4)K_{12}\!(p_1,p_2)=\\
&=\left(\left(K_{34}\star K_{23}\right)\star K_{12}\right)\!(p_1,p_4),
\end{align*}  
which proves associativity of the operation $\star$, whenever it makes sense.

On the other hand, considering a pair $\left(P_1,P_2\right)$ of objects of the category $\mathsf{Bun}_{\calG,M}$, by what was proved in Proposition~\ref{prop-gaugegr} of Subsubsection~\ref{sssec-gaugetrsf}, any bundle automorphism of $P_1$, hence a gauge transformation, corresponds uniquely to an element of $C^{\infty}\!(P_1\odot P_1,\calG)^{\calG^2}$.
In particular, the unique element associated to the identity map on $P_1$ is simply $\phi_{P_1}^{-1}$; namely
\[
K_{\id_{P_1}}(p_1,q_1)=\phi_{P_1}(q_1,p_1)=\phi_{P_1}(p_,q_1)^{-1}.
\]
I want to compute an explicit expression for $K_{12} \star \phi_{P_1}^{-1}$, for any $K_{12}\in C^{\infty}\!(P_1\odot P_2,\calG)^{\calG^2}$:
\begin{align*}
\left(K_{12} \star \phi_{P_1}^{-1}\right)\!(p_1,p_2)&=K_{12}\!(q_1,p_2)\phi_{P_1}\!(p_1,q_1)^{-1}=\\
&\overset{\text{by independence of the choice of $q_1$}}=K_{12}\!(p_1,p_2)\phi_{P_1}\!(p_1,p_1)=\\
&=K_{12}\!(p_1,p_2)\iota_{\calG}\!\left(\varepsilon_1(p_1)\right)=\\
&=K_{12}\!(p_1,p_2),
\end{align*} 
where $\pi_1(p_1)=\pi_1(q_1)$ and since
\[
s_{\calG}\!\left(K_{12}(p_1,p_2)\right)=\epsilon_1(p_1).
\]
On the other hand, using the same notations as before, I compute explicitly $\phi_{P_2}^{-1}\star K_{12}$:
\begin{align*}
\left(\phi_{P_2}^{-1}\star K_{12}\right)\!(p_1,p_2)&=\phi_{P_2}\!(q_2,p_2)^{-1}K_{12}\!(p_1,q_2)=\\
&\overset{\text{by independence of the choice of $q_2$}}=\phi_{P_2}\!(p_2,p_2)^{-1}K_{12}\!(p_1,p_2)=\\
&=\iota_{\calG}\!\left(\varepsilon_2(p_2)\right)K_{12}\!(p_1,p_2)=\\
&=K_{12}\!(p_1,p_2),
\end{align*}
where $\pi_2(p_2)=\pi_2(q_2)$, and by
\[
t_{\calG}\!\left(K_{12}\!(p_1,p_2)\right)=\varepsilon_2(p_2).
\]
Hence, for any object $P_1$ of the category $\mathsf{Bun}_{\calG,M}$, there is an element $\phi_{P_1}^{-1}$, which corresponds to the unit map for the operation $\star$.

At last, for any pair $\left(P_1,P_2\right)$ of objects of the category $\mathsf{Bun}_{\calG,M}$ and any morphism between them represented by $K_{12}\in C^{\infty}\!(P_1\odot P_2,\calG)^{\calG^2}$, Lemma~\ref{lem-gengaugeinv1} of Subsection~\ref{ssec-gengauge} implies that there is another generalized gauge transformation, whose associated bundle morphism is the inverse of the bundle morphism represented by $K_{12}$; I denote this generalized gauge transformation by $\widetilde{K}_{12}$.
Let me compute explicitly the product $\widetilde{K}_{12}\star K_{12}$:
\begin{align*}
\left(\widetilde{K}_{12}\star K_{12}\right)\!(p_1,q_1)&=\widetilde{K}_{12}\!(p_2,q_1)K_{12}\!(p_1,p_2)=\\
&=K_{12}\!(q_1,p_2)^{-1}K_{12}\!(p_1,p_2)=\\
&\overset{\text{by definition of $\phi_{P_1}$}}=K_{12}\!(p_1\phi_{P_1}\!(p_1,q_1),p_2)^{-1}K_{12}\!(p_1,p_2)=\\
&\overset{\text{by $\calG\times \calG$-equivariance of $K_{12}$}}=\phi_{P_1}\!(p_1,q_1)^{-1}K_{12}\!(p_1,p_2)^{-1}K_{12}\!(p_1,p_2)=\\
&=\phi_{P_1}\!(p_1,q_1)^{-1}\iota_{\calG}\!\left(s_{\calG}\!\left(K_{12}\!(p_1,p_2)\right)\right)=\\
&=\phi_{P_1}\!(p_1,q_1)^{-1}\iota_{\calG}\!\left(\varepsilon_1(p_1)\right)=\\
&=\phi_{P_1}\!(p_1,q_1)^{-1},
\end{align*}
where $p_2\in P_2$ is such that $\pi_2(p_2)=\pi_1(p_1)$, and since
\[
s_{\calG}\!\left(\phi_{P_1}\!(p_1,q_1)^{-1}\right)=t_{\calG}\!\left(\phi_{P_1}\!(p_1,q_1)\right)=\varepsilon_(p_1).
\]
On the other hand, similar computations yield
\[
K_{12}\star \widetilde{K}_{12}=\phi_{P_2},
\]
whence the assignment
\[
K_{12}\in C^{\infty}\!(P_1\odot P_2,\calG)^{\calG^2}\leadsto C^{\infty}\!(P_2\odot P_1,\calG)^{\calG^2}\ni \widetilde{K}_{12}
\]
gives an inverse for the operation $\star$.

Putting all these computations together, one sees that the category $\mathsf{Bun}_{\calG,M}$ with principal $\calG$-bundles over $M$ as objects can be made into the set of objects of an abstract groupoid.
Namely, to any pair of objects $\left(P_1,P_2\right)$ of $\mathsf{Bun}_{\calG,M}$, I associate the set
\[
\left(P_1,P_2\right)\leadsto C^{\infty}\!(P_1\odot P_2,\calG)^{\calG^2}
\]
of generalized gauge transformations between $P_1$ and $P_2$.
There are maps $s$, $t$ (the source and target map respectively) from the set of all sets of the form $C^{\infty}\!(P_1\odot P_2,\calG)^{\calG^2}$, for any two objects $P_1$, $P_2$ of $\mathsf{Bun}_{\calG,M}$, to the objects of $\mathsf{Bun}_{\calG,M}$; $\iota$, the unit map, from the objects of the category $\mathsf{Bun}_{\calG,M}$, to the set $C^{\infty}\!(P\odot P,\calG)^{\calG^2}$, for some object $P$ of $\mathsf{Bun}_{\calG,M}$, defined respectively via
\begin{align*}
s\left(K_{12}\right)&\colon=P_1,\quad K_{12}\in C^{\infty}\!(P_1\odot P_2,\calG)^{\calG^2},\\ 
t\left(K_{12}\right)&\colon=P_2,\quad K_{12}\in C^{\infty}\!(P_1\odot P_2,\calG)^{\calG^2},\\ 
i\left(P\right)&\colon=\phi_P^{-1}\in C^{\infty}\!(P\odot P,\calG)^{\calG^2}\cong C^{\infty}\!(P,\calG)^{\calG} .
\end{align*}
There is a partially defined, associative product of the set of sets of the form $C^{\infty}\!(P_1\odot P_2,\calG)^{\calG^2}$:
\begin{align*}
\star :C^{\infty}\!(P_2\odot P_3,\calG)^{\calG^2}\times C^{\infty}\!(P_1\odot P_2,\calG)^{\calG^2}&\to C^{\infty}\!(P_1\odot P_3,\calG)^{\calG^2}\\
\left(K_{23},K_{12}\right)&\mapsto K_{23}\star K_{12}.
\end{align*}
It is obvious that 
\begin{equation}\label{eq-sourtarid}
\begin{aligned}
s\left(K_{23}\star K_{12}\right)&=P_1=s\left(K_{12}\right),\\
t\left(K_{23}\star K_{12}\right)&=P_3=t\left(K_{23}\right),\\
s\left(i\left(P\right)\right)&=P=t\left(i\left(P\right)\right);\\
i\left(t\left(K_{12}\right)\right)\star K_{12}&=K_{12},\quad K_{12}\star i\left(s\left(K_{12}\right)\right)=K_{12},\quad \forall K_{12}\in C^{\infty}\!\left(P_1\odot P_2,\calG\right)^{\calG^2}.
\end{aligned}
\end{equation}
It was also proved that there exists, for any $K_{12}\in C^{\infty}\!(P_1\odot P_2,\calG)^{\calG^2}$, a unique element, which was previously denoted by $\widetilde{K}_{12}\in C^{\infty}\!(P_2\odot P_1,\calG)^{\calG^2}$, which satisfies the property
\[
K_{12}\star \widetilde{K}_{12}=\phi_{P_2}=i\left(t\left(K_{12}\right)\right),\quad \widetilde{K}_{12}\star K_{12}=\phi_{P_1}=i\left(s\left(K_{12}\right)\right).
\]
Hence, it makes sense to define the inversion map
\[
j\!\left(K_{12}\right)\colon=\widetilde{K}_{12},\quad\forall K_{12}\in C^{\infty}\!\left(P_1\odot P_2,\calG\right)^{\calG^2}.
\]
The {\em groupoid of generalized gauge transformations} is denoted by $C^{\infty,\calG^2}$, so that 
\begin{equation}\label{eq-inversion}
C^{\infty,\calG^2}_{P_1,P_2}\colon=C^{\infty}\!\left(P_1\odot P_2,\calG\right)^{\calG^2},\quad C^{\infty,\calG^2}_{P}\cong C^{\infty}\!\left(P,\calG\right)^\calG,
\end{equation}
the latter being, as was said before, a consequence of Proposition~\ref{prop-gaugegr} of Subsubsection~\ref{sssec-gaugetrsf}; the source, target and unit map of $C^{\infty,\calG^2}$ are defined in (\ref{eq-sourtarid}).

Summarizing all the computations so far, I can finally state the following
\begin{Thm}\label{thm-gaugegroupoid}
The $6$-tuple $\left(C^{\infty,\calG^2},\mathsf{Bun}_{\calG,M},s,t,i,j\right)$, for any manifold $M$ and any groupoid $\calG$, is an abstract groupoid in the sense of Definition~\ref{def-groupoid}, where the source, target and unit map are defined in (\ref{eq-sourtarid}), and the inversion map is defined in (\ref{eq-inversion}); it is obviously isomorphic to the abstract groupoid of bundle morphisms with set of objects $\mathsf{Bun}_{\calG,M}$ by Theorem~\ref{thm-gengaugeeq} of Subsection~\ref{ssec-gengauge}.
For any object $P$ of $\mathsf{Bun}_{\calG,M}$, the isotropy group $C^{\infty,\calG^2}_P$ is isomorphic to the gauge group of $P$, which is denoted $C^{\infty}\!\left(P,\calG\right)^{\calG}$.
\end{Thm}

\section{Equivariant morphisms between Hilsum--Skandalis morphisms and a subgroupoid of the groupoid of generalized gauge transformations}\label{sec-equivhilsumskand}
In this section, I discuss the notion of {\em Hilsum--Skandalis morphism} or {\em generalized morphism} between two groupoids $\calG$ and $\calH$, which are particular right principal $\calH$-bundles with an additional left $\calG$-action compatible with projection, momentum and right $\calH$-action.
There is a natural notion of morphisms between such bundles, which is a specialization of the concept of bundle morphisms in the sense of Definition~\ref{def-equivbun}; therefore, it is possible to introduce also a notion of generalized gauge transformations of Hilsum--Skandalis morphisms, and I will characterize completely such generalized gauge transformations, proving that they form a subgroupoid of the groupoid of generalized gauge transformations of right principal $\calH$-bundles, introduced and discussed at the end of the previous section.

\subsection{The notion of Hilsum--Skandalis morphisms and their division map}\label{ssec-hilskandmor}
In this Subsection, I define Hilsum--Skandalis morphisms between two groupoids $\calG$ and $\calH$; as we will see, such morphisms are in fact right principal $\calH$-bundles with a compatible $\calG$-action, and their name seems somehow mysterious, because there is apparently no notion of morphism in the strict sense.
In truth, hidden in the local structure of a Hilsum--Skandalis morphism as a bundle, there is a ``stretched'' notion of morphism.
This is particularly evident in~\cite{Con},~\cite{HS}, where such bundles were introduced using a {\em local} point of view; in~\cite{C2}, I discuss in great detail the local properties of Hilsum--Skandalis morphisms (although I prefer to use the denomination ``generalized morphisms'' in~\cite{C2}), and it is clarified there that a local generalized morphism from a groupoid $\calG$ to a groupoid $\calH$ (which has, to a certain extent, the properties of a true morphism) encodes all the data one needs to provide a left $\calG$-action on a previously constructed (from local data) right principal $\calH$-bundle.
In fact, a morphism between groupoids $\calG$ and $\calH$ in the sense of Definition~\ref{def-morgroupoid} gives rise in a natural way to a generalized morphism in the sense that I give now.
\begin{Def}\label{def-hilsumskan}
Given two Lie groupoids $\calG$ and $\calH$, a {\em Hilsum--Skandalis morphism from $\calG$ to $\calH$} (from now on, shortly HS morphism) is a right principal $\calH$-bundle $\left(P,\pi,\varepsilon,X_{\calG}\right)$ over $X_{\calG}$, obeying the following additional conditions:
\begin{itemize}
\item[i)] the pair $\left(P,\pi\right)$ defines a left $\calG$-action on $P$ with momentum $\pi$;
\item[ii)] the momentum $\varepsilon$ for the right $\calH$-action is $\calG$-invariant, and both actions are compatible in the sense that
\[
(gp)h=g(ph),\quad s_{\calG}(p)=\pi(p),\quad t_{\calH}(h)=\varepsilon(p).
\]
\end{itemize}
(Notice that the $\calG$-invariance of the momentum $\varepsilon$ makes both sides of the compatibility condition between both actions well-defined.)
\end{Def}

\begin{Exa}
(For a detailed discussion of the first two examples of HS morphisms, I refer to~\cite{C2})
\begin{itemize}
\item[i)] Consider two Lie groups $G$ and $H$ as Lie groupoids, with trivial unit, source and target maps; then a HS morphism from $G$ to $H$ is a morphism of Lie groups in the usual sense. 
\item[ii)] Considering again two Lie groups $G$ and $H$, as in the previous example, acting this time respectively (from the left) on the manifolds $M$ and $N$, one can consider the associated action groupoids $G\ltimes M$ and $H\ltimes N$.
Then, HS morphisms from $G\ltimes M$ to $H\ltimes N$ are in one-to-one correspondence with right $H$-principal bundles $P$ over $M$ endowed with a lift of the left $G$-action on $M$ to $P$ and to the associated bundle $P\times_H N$ and with a $G$-equivariant global section of $P\times_H N$.
\item[iii)] (See e.g.~\cite{L-GTX} for more details) Given a Lie group $G$, viewed as a (trivial) Lie groupoid, and a Lie groupoid $\Gamma$ over the manifold of objects $\Gamma_0$, HS morphisms from $\Gamma$ to $G$ are called {\em principal $G$-bundles over the groupoid $\Gamma$}.
Clearly, such HS morphisms, for $\Gamma=H\ltimes M$, $H$ being a Lie group operating from the left on the manifold $M$ (which is $\Gamma_0$), correspond to $H$-equivariant principal $G$-bundles over $M$. 
\end{itemize}
\end{Exa}

In particular, it follows from Definition~\ref{def-hilsumskan} that a HS morphism $P$ from $\calG$ to $\calH$ possesses a division map $\phi_P$ in the sense explained at the beginning of Section~\ref{sec-princbun}.
Clearly, there is more at work in the presence of HS morphisms, and it is natural to expect that the division map of a HS morphism has particular features.
The first important fact in this direction is encoded in the following
\begin{Lem}\label{lem-prodhilskand}
Given two groupoids $\calG$ and $\calH$, and two HS morphisms $P$ and $\widetilde{P}$ from $\calG$ to $\calH$ in the sense of Definition~\ref{def-hilsumskan}, the product bundle $\left(P\times\widetilde{P},\pi\times \widetilde{\pi},\varepsilon\times\widetilde{\varepsilon},X_{\calG}\times X_{\calG}\right)$, where I have borrowed the notations from Subsection~\ref{ssec-fibprodbun}, is a HS morphism from $\calG^2$ to $\calH^2$. 
\end{Lem} 
Notice that Corollary~\ref{cor-prodbun} implies already that $\left(P\times\widetilde{P},\pi\times \widetilde{\pi},\varepsilon\times\widetilde{\varepsilon},X_{\calG}\times X_{\calG}\right)$ is a right principal $\calH^2$-bundle over $X_{\calG}\times X_{\calG}$.
Lemma~\ref{lem-prodhilskand} is a trivial consequence of the more general statement contained in 
\begin{Lem}\label{lem-prodhilskand1}
Let $\calG_1$, $\calG_2$, $\calH_1$ and $\calH_2$ four groupoids; let $P_1$, resp.\ $P_2$, a HS morphism in the sense of Definition~\ref{def-hilsumskan} from $\calG_1$ to $\calH_1$, resp.\ $\calG_2$ to $\calH_2$, with projection and momentum labelled by $\pi_1$, $\varepsilon_1$ and $\pi_2$, $\varepsilon_2$ respectively.

Then, the $4$-tuple $\left(P_1\times P_2,\pi_1\times\pi_2,\varepsilon_1\times\varepsilon_2,X_{\calG_1}\times X_{\calG_2}\right)$ is a HS morphism from $\calG_1\times\calG_2$ to $\calH_1\times\calH_2$, the product of the HS morphisms $P_1$ and $P_2$.
\end{Lem}
\begin{proof}
Lemma~\ref{lem-prodbun} shows already that the above $4$-tuple is a right principal $\calH_1\times\calH_2$-bundle, with projection $\pi_1\times \pi_2$ and momentum $\varepsilon_1\times\varepsilon_2$; it remains to show that there is a left $\calG_1\times\calG_2$-action on $P_1\times P_2$ with momentum $\pi_1\times\pi_2$.
In fact, consider a pair $(g_1,g_2)$ in the product groupoid $\calG_1\times\calG_2$ and a pair in $P_1\times P_2$, such that
\[
\pi_1(p_1)=s_{\calG_1}(g_1),\quad \pi_2(p_2)=s_{\calG_2}(g_2);
\]
then, a natural left $\calG_1\times\calG_2$-action is
\[
(g_1,g_2)(p_1,p_2)\colon=(g_1p_1,g_2p_2).
\]
The above action is clearly well-defined, and moreover
\begin{align*}
t_{\calG_1\times\calG_2}(g_1,g_2)&=(t_{\calG_1}(g_1),t_{\calG_2}(g_2))=\\
&=\left(\pi_1(g_1p_1),\pi_2(g_2p_2)\right)=\\
&=\left(\pi_1\times\pi_2\right)(g_1p_1,g_2p_2).
\end{align*}
On the other hand, it holds
\begin{align*}
\left((g_1,g_2)(p_1,p_2)\right)(h_1,h_2)&=(g_1p_1,g_2p_2)(h_1,h_2)=\\
&=\left((g_1p_1)h_1,(g_2p_2)h_2\right)=\\
&=\left(g_1(p_1h_1),g_2(p_2h_2)\right)=\\
&=(g_1,g_2)\left((p_1,p_2)(h_1,h_2)\right),\\
&s_{\calG_1}(g_1)=\pi_1(p_1),\quad t_{\calH_1}(h_1)=\varepsilon_1(p_1),\\
&s_{\calG_2}(g_2)=\pi_2(p_2),\quad t_{\calH_2}(h_2)=\varepsilon_2(p_2),
\end{align*}
which shows that both actions are compatible.
It remains to show that the momentum $\varepsilon_1\times\varepsilon_2$ is left $\calG_1\times\calG_2$-invariant:
\begin{align*}
(\varepsilon_1\times\varepsilon_2)(g_1p_1,g_2p_2)&=\left(\varepsilon_1(g_1p_1),\varepsilon_2(g_2p_2)\right)=\\
&=\left(\varepsilon_1(g_1),\varepsilon_2(g_2)\right)=\\
&=(\varepsilon_1\times\varepsilon_2)(p_1,p_2),\quad s_{\calG_i}(g_i)=\pi_i(p_i),\quad i=1,2.
\end{align*}
\end{proof}
Analogously to what I did already in Subsection~\ref{ssec-fibprodbun}, I consider the special case $\calG_1=\calG_2=\calG$ and $\calH_1=\calH_2=\calH$, and I consider two HS morphisms $P$ and $\widetilde{P}$ from $\calG$ to $\calH$; I further consider the diagonal $\Delta_{\calG}\colon=\Delta_{X_{\calG}}$, and I use the same notation for the diagonal imbedding of the diagonal into $X_{\calG}\times X_{\calG}$.
I consider then the restriction of the HS morphism $P\times\widetilde{P}$ from $\calG^2$ to $\calH^2$ to the diagonal $\Delta_{\calG}$; Lemma~\ref{lem-fibprod} implies that the fibred product $P\odot\widetilde{P}$ is a right principal $\calH^2$-bundle, with projection 
\[
\overline{\pi}(p,\widetilde{p})=\pi(p)=\widetilde{\pi}(\widetilde{p}),
\]
and momentum
\[
\overline{\varepsilon}(p,\widetilde{p})\colon=\left(\varepsilon{p},\widetilde{\varepsilon}(\widetilde{p})\right).
\]
As one could expect quite naturally, the following holds: 
\begin{Lem}\label{lem-fibprodhils}
Given two groupoids $\calG$ and $\calH$, the fibred product bundle $P\odot\widetilde{P}$ of two HS morphisms $P$ and $\widetilde{P}$ in the sense of Definition~\ref{def-hilsumskan} from $\calG$ to $\calH$ is a HS morphism from $\calG$ to $\calH^2$.
\end{Lem}
\begin{proof}
One has only to show that there is a left $\calG$-action on $P\odot\widetilde{P}$ with momentum $\overline{\pi}$, which is compatible with the right $\calH^2$-action and such that the momentum $\overline{\varepsilon}$ of the right $\calH^2$-action is $\calG$-invariant.
Consider a pair $(p,\tildep)$ in $P\odot \widetilde{P}$ and an element $g\in\calG$, such that
\[
\pi(p)=\widetilde{\pi}(\tildep)=s_{\calG}(g); 
\]
then, define the left $\calG$-action as the restriction to the diagonal of $\calG^2$ of the left $\calG^2$-action on $P\times \widetilde{P}$:
\[
g(p,\tildep)\colon=(gp,g\tildep).
\]
It is immediate to verify that the above formula defines a left $\calG$-action on $P\odot\widetilde{P}$; a slight modification of the arguments used in the proof of Lemma~\ref{lem-prodhilskand1} shows that this action is compatible with the right $\calH^2$-action and that $\overline{\varepsilon}$ is $\calG$-invariant.
\end{proof}

Consider now a HS morphism $P$ from $\calG$ to $\calH$; since $P$ is a right principal $\calH$-bundle, it possesses a (uniquely defined) division map $\phi_P$.
In the following proposition are listed all properties of $\phi_P$:
\begin{Prop}\label{prop-proddivhils}
Given a HS morphism $P$ from the groupoid $\calG$ to the groupoid $\calH$, the division map $\phi_P$ is a map from the fibred product bundle $P\odot P$ to $\calH$ with the following additional properties:
\begin{itemize}
\item[i)] (Compatibility with momentum)
\[
t_{\calH}\circ \phi_P=\varepsilon\circ \pr_1,\quad s_{\calH}\circ \phi_P=\varepsilon\circ \pr_2.
\]
\item[ii)] On the diagonal submanifold of the total space of $P\odot P$ holds
\[
\phi_P(p,p)=\iota_{\calH}(\varepsilon(p)),\quad \forall p\in P.
\] 
\item[iii)] for any pair $(p,q)\in P\odot P$, the following equation holds
\[
\phi_P(p,q)=\phi_P(q,p)^{-1};
\] 
notice that the previous equation makes sense, since $(p,q)\in P\odot P$ implies that $(q,p)\in P\odot P$ also.
\item[iv)] The division map $\phi_P$ is a $\calG$-invariant, $\calH^2$-equivariant map from the fibred product bundle $P\odot P$ to $\calH$, endowed with the right generalized conjugation defined in Remark~\ref{rem-conjright} by the pair of maps $\left(J_{\conj}^R,\Psi_{\conj}^R\right)$ in Subsubsection~\ref{sssec-genconjgroupoid}.  
\end{itemize}
\end{Prop}
\begin{proof}
Almost all properties of the division map $\phi_P$ follows from Proposition~\ref{prop-prophi}; it remains to show that $\phi_P$ is $\calG$-invariant, i.e.\ one has to show
\[
\phi_P(gp_1,gp_2)=\phi_P(p_1,p_2),\quad \forall (p_1,p_2)\in P\odot P,\quad g\in\calG,\quad s_{\calG}(g)=\pi(p_1)=\pi(p_2).
\]
(Notice that the $\calG$-invariance of the momentum $\varepsilon$ makes the above identity well-defined.)
In fact, by its very definition, the division map $\phi_P$ satisfies the identity
\[
gp_2=gp_1\phi_P(gp_1,gp_2)=g\left(p_1\phi_P(p_1,p_2)\right),
\]
using the compatibility of both actions.
Since the right $\calH$-action is free, the claim follows.
\end{proof}

\subsection{Morphisms between HS morphisms}\label{ssec-morhilskan}
As I have already pointed out, a HS morphism $P$ from a morphism $\calG$ to $\calH$ is in particular a right principal $\calH$-bundle.
Given now two HS morphisms $P_1$ and $P_2$ from $\calG$ to $\calH$, it is therefore natural to consider morphisms from $P_1$ to $P_2$ in the sense of Definition~\ref{def-equivbun} in Subsection~\ref{ssec-eqmorph}, namely fibre-preserving, $\calH$-equivariant maps from $P_1$ to $P_2$.
Since $P_1$ and $P_2$ are both also left $\calG$-spaces, and such a morphism is immediately momentum-preserving (w.r.t.\ the momenta $\varepsilon_1$ and $\varepsilon_2$ of the right $\calH$-actions on $P_1$ and $P_2$ respectively), it makes to consider the following subset of the morphisms in the sense of Definition~\ref{def-equivbun}:
\begin{Def}\label{def-equivhils}
Given two Lie groupoids $\calG$ and $\calH$, any two HS morphisms $P_1$ and $P_2$ from $\calG$ to $\calH$ in the sense of Definition~\ref{def-hilsumskan}, a morphism $\sigma$ from $P_1$ to $P_2$ is a morphism from the right $\calH$-bundle $P_1$ to the right $\calH$-bundle $P_2$ in the sense of Definition~\ref{def-equivbun}, which satisfies additionally the following requirement ($\calG$-equivariance):
\[
\sigma(gp_1)=g\sigma(p_1),\quad \forall p_1\in P_1,\quad g\in\calG,\quad s_{\calG})g_=\pi_1(p_1).
\]
\end{Def} 
It is immediate to verify that, given three HS morphisms $P_1$, $P_2$ and $P_3$ from $\calG$ to $\calH$, and morphisms $\sigma_{12}$ from $P_1$ to $P_2$ and $\sigma_{23}$ from $P_2$ to $P_3$ in the sense of Definition~\ref{def-equivhils}, the composition $\sigma_{23}\circ\sigma_{12}$ is again a morphism of HS morphisms from $P_1$ to $P_3$ in the sense of Definition~\ref{def-equivhils}.
It is also easy to verify that the identity map $\id_P$ of a HS morphism is a morphism of HS morphisms.
Therefore, given two groupoids $\calG$ and $\calH$, analogously to what I did in Subsection~\ref{ssec-eqmorph}, I consider the category $\mathsf{HS}_{\calG,\calH}$, whose ingredients are as follows:
\begin{itemize}
\item[i)] {\bf Objects}: objects of the category $\mathsf{HS}_{\calG,\calH}$ are HS morphisms from $\calG$ to $\calH$;
\item[ii)] {\bf Morphisms}: morphisms of the category $\mathsf{HS}_{\calG,\calH}$ are morphisms of HS morphisms in the sense of Definition~\ref{def-equivhils}.
\end{itemize}  
It is clear that the category $\mathsf{HS}_{\calG,\calH}$ is a subcategory of $\mathsf{Bun}_{X_{\calG},\calH}$.
As a consequence, Lemma~\ref{lem-inverequiv} shows that every morphism in $\mathsf{HS}_{\calG,\calH}$ is bijective; moreover, the same machinery developed in Subsection~\ref{ssec-gengauge} can be applied to morphisms in $\mathsf{HS}_{\calG,\calH}$ to show that every such morphism is invertible in the same category.
This is what I am going to do in what follows.

Consider two objects $P_1$, $P_2$ in $\mathsf{HS}_{\calG,\calH}$ and a morphism $\sigma$ of the category $\mathsf{HS}_{\calG,\calH}$ between them; viewing $\mathsf{HS}_{\calG,\calH}$ as a subcategory of $\mathsf{Bun}_{X_{\calG},\calH}$, to $\sigma$ belongs a unique generalized gauge transformation $K_{\sigma}$ in $C^{\infty}\!\left(P_1\odot P_2,\calH\right)^{\calH^2}$, defined via
\[
K_{\sigma}(p_1,p_2)=\phi_{P_2}\!\left(p_2,\sigma(p_1)\right)=\phi_{P_1}\!\left(\sigma^{-1}(p_2),p_1\right),
\]
where the second equality is a consequence of Theorem~\ref{thm-gaugeinvdiv}.
Let me point out a caveat: the inverse $\sigma^{-1}$ in the above identity is {\em not} the inverse of $\sigma$ in the category $\mathsf{HS}_{\calG,\calH}$, but the inverse of $\sigma$ in $\mathsf{Bun}_{X_{\calG},\calH}$.
For more details about generalized gauge transformations, I refer to Subsection~\ref{ssec-gengauge}.
Let me introduce a new notion at this point
\begin{Def}\label{def-gengaugehils}
Given two groupoids $\calG$ and $\calH$, and two objects $P_1$, $P_2$ in $\mathsf{HS}_{\calG,\calH}$, a {\em Hilsum--Skandalis generalized gauge transformation} (shortly, a HS generalized gauge transformation) between $P_1$ and $P_2$ is a $\calG$-invariant, $\calH^2$-equivariant map from the fibred product bundle $P_1\odot P_2$ to $\calH$; as in Definition~\ref{def-gengauge}, $\calH$ is a right $\calH^2$-space w.r.t.\ the right generalized conjugation introduced and discussed in Subsubsection~\ref{sssec-genconjgroupoid}, while $P_1\odot P_2$ is a left $\calG$-space and a right $\calH^2$-bundle in virtue of Lemma~\ref{lem-fibprodhils}.
The set of HS generalized gauge transformations between the HS morphisms $P_1$ and $P_2$ from $\calG$ to $\calH$ is denoted by $C^{\infty}_\calG\!\left(P_1\odot P_2,\calH\right)^{\calH^2}$.
\end{Def} 
To see more explicitly the properties of a HS generalized gauge transformation between $P_1$ and $P_2$ in the objects of $\mathsf{HS}_{\calG,\calH}$, I refer to Remark~\ref{rem-explgengauge} with a caveat: in the case of a HS generalized gauge transformation $K$, one has to consider the additional equation
\begin{equation}
K(gp_1,gp_2)=K(p_1,p_2),\quad \forall (p_1,p_2)\in P_1\odot P_2,\quad g\in\calG,\quad s_{\calG}(g)=\pi_1(p_1)=\pi_2(p_2),
\end{equation}
namely the $\calG$-invariance of $K$.

Theorem~\ref{thm-gengaugeeq} has a natural counterpart in the framework of HS generalized gauge transformations, namely
\begin{Thm}\label{thm-gengaugehils}
Given two groupoids $\calG$ and $\calH$, the set of morphisms $\mathsf{Mor}_{\mathsf{HS}_{\calG,\calH}}\!(P_1,P_2)$ from the object $P_1$ to the object $P_2$ in $\mathsf{Ob}\!\left(\mathsf{HS}_{\calG,\calH}\right)$ is in one-to-one correspondence with the set $C^{\infty}_\calG\!\left(P_1\odot P_2,\calH\right)^{\calH^2}$ of HS generalized gauge transformations between $P_1$ and $P_2$. 
\end{Thm}
\begin{proof}
Since a morphism $\sigma$ from the HS morphism $P_1$ to the HS morphism $P_2$ is in particular a morphism from the right $\calH$-bundle $P_1$ to the right $\calH$-bundle $P_2$ in the sense of Definition~\ref{def-equivbun}, it follows immediately by Theorem~\ref{thm-gengaugeeq} that the assignment
\[
\sigma\leadsto K_\sigma(p_1,p_2)\colon=\phi_{P_2}(p_2,\sigma(p_1)),\quad (p_1,p_2)\in P_1\odot P_2 
\]
is a generalized gauge transformation between $P_1$ and $P_2$ in the sense of Definition~\ref{def-gengauge}, i.e.\ an element of $C^{\infty}\!\left(P_1\odot P_2,\calH\right)^{\calH^2}$.
To prove that $K_\sigma$ is a HS generalized gauge transformation, it suffices to prove that it is $\calG$-invariant.
This follows in turn from the $\calG$-invariance of $\sigma$ and from Point $iv)$ of Proposition~\ref{prop-proddivhils}.

On the other hand, consider now a HS generalized gauge transformation $K$ between $P_1$ and $P_2$ in the sense of Definition~\ref{def-gengaugehils}; $K$ is in particular a generalized gauge transformation between the right $\calH$ bundles $P_1$ and $P_2$ in the sense of Definition~\ref{def-gengauge}.
Thus, by Theorem~\ref{thm-gengaugeeq}, the following assignment defines a morphism from the right $\calH$-bundle $P_1$ to the right $\calH$-bundle $P_2$:
\[
\sigma_K(p_1)\colon=p_2K(p_1,p_2),\quad \pi_1(p_1)=\pi_2(p_2).
\]
Recall that the definition of $\sigma_K$ does not depend on the choice of $p_2$ such that $\pi_1(p_1)=\pi_2(p_2)$.
To show that $\sigma_K$ is a morphism of HS morphisms, it remains to show that it is $\calG$-invariant, i.e.\ one has to show 
\[
\sigma_K(gp_1)=g\sigma_K(p_1),\quad \forall p_1\in P_1,\quad g\in\calG,\quad s_{\calG}(g)=\pi_1(p_1).
\]
But this is a consequence of the following arguments: by definition,
\begin{align*}
\sigma_K(gp_1)=\tildep_2K(gp_1,\tildep_2),\quad \pi_2(\tildep_2)=\pi_1(gp_1)=t_{\calG}(g).
\end{align*}
Since the previous formula does not depend on the choice of the representative $\tildep_2$, one can choose 
\[
\tildep_2=g p_2,\quad \pi_1(p_1)=\pi_2(p_2),
\]
and from this it follows, by $\calG$-invariance of $K$,
\begin{align*}
\sigma_K(gp_1)&=\tildep_2K(gp_1,\tildep_2)=\\
&=gp_2K(gp_1,gp_2)=\\
&=gp_2K(p_1,p_2)=\\
&=g\sigma_K(p_1),
\end{align*} 
whence the claim follows.

Obviously, since the assignments of a HS generalized gauge transformation $K$ to a morphism $\sigma_K$ between HS morphisms and viceversa are constructed by the same rules as in the proof of Theorem~\ref{thm-gengaugeeq}, it follows immediately by the very same arguments that 
\[
\mathsf{Mor}_{\mathsf{HS}_{\calG,\calH}}\!(P_1,P_2)\ni\sigma\leadsto K_\sigma\in C^{\infty}_\calG\!\left(P_1\odot P_2,\calH\right)^{\calH^2}
\]
and
\[
C^{\infty}_\calG\!\left(P_1\odot P_2,\calH\right)^{\calH^2}\ni K\leadsto \sigma_K\in\mathsf{Mor}_{\mathsf{HS}_{\calG,\calH}}\!(P_1,P_2)
\]
are inverse to each other, hence proving the Theorem.
\end{proof}
At this point, we know that $i)$ a morphism $\sigma$ in the category $\mathsf{HS}_{\calG,\calH}$ is bijective and $ii)$ any such morphism corresponds uniquely to a HS generalized gauge transformation $K_\sigma$.
Viewing such a HS generalized gauge transformation $K_\sigma$ as a generalized gauge transformation, Lemma~\ref{lem-gengaugeinv1} implies that
\[
K_{\sigma^{-1}}\!\left(p_2,p_1\right)\colon=K_{\sigma}(p_1,p_2)^{-1}
\]
is a generalized gauge transformation between $P_2$ and $P_1$, viewed both as objects of $\mathsf{Bun}_{\calH,M}$.
The $\calG$-invariance of $K_{\sigma}$ implies immediately that $K_{\sigma^{-1}}$ is also $\calG$-invariant; thus, $K_{\sigma^{-1}}$ is a HS generalized gauge transformation.
Lemma~\ref{lem-gengaugeinv2} implies immediately that the morphism $\sigma_{K_{\sigma^{-1}}}$ in $\mathsf{HS}_{\calG,\calH}$ corresponding to the HS generalized gauge transformation $K_{\sigma^{-1}}$ between $P_2$ by $P_1$ by Theorem~\ref{thm-gengaugehils} is the inverse of $\sigma$, thus proving that every morphism in $\mathsf{HS}_{\calG,\calH}$ is invertible.

Let me consider now an object $P$ in $\mathsf{HS}_{\calG,\calH}$ and a morphism $\sigma$ in $\mathsf{Mor}_{\mathsf{HS}_{\calG,\calH}}\!(P,P)$: Theorem~\ref{thm-gengaugehils} ensures the existence of a HS generalized gauge transformation $K_\sigma$ on $P$.
On the other hand, since in particular $\sigma$ preserves the projection $\pi$, there is a unique map $G_\sigma$ from $P$ to $\calG$, such that
\[
\sigma(p)=pG_\sigma(p),
\] 
and this map can obviously be written in the following form, by the properties of the division map of $P$
\[
G_\sigma(p)=\phi_P(p,\sigma(p))\overset{!}=K_\sigma(p,p),
\]
by the very construction of $K_\sigma$ according to Theorem~\ref{thm-gengaugehils}, or, in other words, $G_\sigma$ is the restriction to the diagonal in $P\odot P$.
By the very properties of a HS generalized gauge transformation, it follows immediately that the map $G_\sigma$ enjoys the following properties:
\begin{itemize}
\item[i)] (Compatibility with momentum)
\[
t_{\calH}\circ G_{\sigma}=s_{\calG}\circ G_\sigma=\varepsilon.
\] 
\item[ii)] ($\calG$-invariance and $\calH$-equivariance)
\begin{align*}
G_\sigma(gp)&=G_\sigma(p),\quad \forall p\in P,\quad\forall g\in \calG,\quad s_{\calG}(g)=\pi(p);\\
G_\sigma(ph)&=h^{-1}G_{\sigma}(p)h,\quad\forall p\in P,\quad\forall h\in\calH,\quad t_{\calH}(h)=\varepsilon(p).
\end{align*}
\end{itemize}
On the other hand, given an object $P$ in $\mathsf{HS}_{\calG,\calH}$, if we consider a HS gauge transformation $G$ of $P$, i.e.\ a map $G$ from $P$ to $\calH$ satisfying
\begin{align*}
t_{\calG}\circ G&=s_{\calG}\circ G=\varepsilon,\\
G(gp)&=G(p),\quad \forall p\in P,\quad\forall g\in \calG,\quad s_{\calG}(g)=\pi(p)\quad\text{and}\\
G(ph)&=h^{-1}G(p)h,\quad\forall p\in P,\quad\forall h\in\calH,\quad t_{\calH}(h)=\varepsilon(p),
\end{align*}
the assignment
\[
G\leadsto \sigma_G(p)\colon=pG(p),
\]
which is well-defined by the first requirement $G$ satisfies, defines in an obvious way a morphism in $\mathsf{HS}_{\calG,\calH}$.
I denote by $C^{\infty}_\calG\!\left(P,\calH\right)^{\calH}$ the set of HS gauge transformation of $P$ in $\mathsf{Ob}\!\left(\mathsf{HS}_{\calG,\calH}\right)$.
By Theorem~\ref{thm-gengaugehils}, there is a (uniquely defined) HS generalized gauge transformation $K_G=K_{\sigma_G}$ on $P$, defined by
\begin{align*}
K_G(p,\tildep)&=\phi_P(p,\sigma_G(\tildep))=\\
&=\phi_P\!\left(p,\tildep G(\tildep)\right)=\\
&=\phi_P(p,\tildep)G(\tildep),\quad\forall (p,\tildep)\in P\odot P.
\end{align*}
It is immediate to verify that the assignments
\[
C^{\infty}_\calG\!\left(P,\calH\right)^{\calH}\ni G\leadsto K_G\in C^{\infty}_\calG\!\left(P\odot P,\calH\right)^{\calH^2}
\]
and
\[
C^{\infty}_\calG\!\left(P\odot P,\calH\right)^{\calH^2}\ni\leadsto \Delta_P^*K\in C^{\infty}_\calG\!\left(P,\calH\right)^{\calH},
\]
where $\Delta_P$ denotes the imbedding of the diagonal of $P\odot P$, are inverse to each other.
Moreover, it is immediate to verify that the pointwise product of two HS gauge transformations of $P$ is again a HS gauge transformation of $P$, and obviously the map
\[
\iota_{\calH}\circ \varepsilon
\]
is also a HS gauge transformation of $P$; the $\calG$-invariance follows from the $\calG$-invariance of the momentum, and the pointwise left and right multiplication by $\iota_{\calH}\circ\varepsilon$ of any HS gauge transformation $G$ fixes $G$.
Finally, the pointwise inverse $G^{-1}$ of a HS generalized gauge transformation $G$ of $P$ is again a HS generalized gauge transformation, and the pointwise product of $G$ with its inverse $G^{-1}$ equals $\iota_{\calH}\circ\varepsilon$.
(The computations so far mimic those of Subsubsection~\ref{sssec-gaugetrsf}.)
Therefore, Proposition~\ref{prop-gaugegr} has the following analogon in the framework of HS morphisms:
\begin{Prop}\label{prop-gaugegrhils}
For any object $P$ in $\mathsf{HS}_{\calG,\calH}$, the set $C^{\infty}_{\calG}\!(P,\calH)^\calH$ of HS gauge transformations from $P$ to $\calG$ is in one-to-one correspondence via the maps
\[
C^{\infty}_{\calG}\!(P,\calH)^\calH\ni G\leadsto K_G\in C^{\infty}_{\calG}(P\odot P,\calH)^{\calH^2}\ni K\leadsto \Delta_P^*\!K\in C^{\infty}_{\calG}(P,\calH)^\calH
\]
with the set of HS generalized gauge transformations $C^{\infty}_{\calG}(P\odot P,\calH)^{\calH^2}$.
Moreover, the set $C^\infty_{\calG}!(P,\calH)^\calH$ is a group, called the {\em HS gauge group of $P$}; thus, the $\mathsf{Mor}_{\mathsf{HS}_{\calG,\calH}}\!(P,P)$, being in one-to-one correspondence with the HS gauge group $C^\infty_{\calG}(P,\calH)^\calH$, inherits a group structure via composition, and the map $\sigma\mapsto G_\sigma$, for any morphism $\sigma$ of $P$, is an isomorphism of groups.
\end{Prop}

I consider now three objects $P_1$, $P_2$ and $P_3$ in $\mathsf{HS}_{\calG,\calH}$, and a HS generalized gauge transformation $K_{12}$ between $P_1$ and $P_2$ and a HS generalized gauge transformation $K_{23}$ between $P_2$ and $P_3$; since they are generalized gauge transformation between $P_1$ and $P_2$ and $P_2$ and $P_3$ respectively, in the sense of Definition~\ref{def-gengauge}, one can consider their product $\star$, as defined at the beginning of Subsection~\ref{ssec-gengaugegroupoid}:
\[
\left(K_{23}\star K_{12}\right)\!(p_1,p_3)\colon=K_{23}(p_2,p_3)K_{12}(p_1,p_2),
\] 
where $\pi_1(p_1)=\pi_2(p_2)=\pi_3(p_3)$, and the result is a generalized gauge transformation between $P_1$ and $P_3$.
I claim now that $K_{23}\star K_{12}$ is $\calG$-invariant, hence the product $\star$ on generalized gauge transformations between right $\calH$-bundles on $X_{\calG}$ descends to a product on HS generalized gauge transformations between objects of $\mathsf{HS}_{\calG,\calH}$.
In fact, the key point is that the above formula does not depend on the choice of $p_2$ in $P_2$, such that $\pi_1(p_1)=\pi_2(p_2)=\pi_3(p_3)$.
Thus, consider $(p_1,p_3)$ in the fibred product bundle $P_1\odot P_3$, and $g\in\calG$, such that
\[
s_{\calG}(g)=\pi_1(p_1)=\pi_3(p_3);
\]
then, it holds 
\[
\pi_1(gp_1)=\pi_3(gp_3)=t_{\calG}(g)=\pi_2(gp_2),
\]
for any $p_2\in P_2$, such that $\pi_2(p_2)=s_{calG}(g)=\pi_1(p_1)=\pi_3(p_3)$.
Therefore, it holds:
\begin{align*}
\left(K_{23}\star K_{12}\right)\!(gp_1,gp_3)&=K_{23}(\tildep_2,gp_3)K_{12}(gp_1,\tildep_2)=\quad \pi_2(\tildep_2)=t_{\calG}(g)\\
&=K_{23}(gp_2,gp_3)K_{12}(gp_1,gp_2)=\\
&=K_{23}(p_2,p_3)K_{12}(p_1,p_2)=\\
&=\left(K_{23}\star K_{12}\right)\!(p_1,p_3),
\end{align*}
and the second equality follows by the independence of the choice of $\tildep_2$.
The product $\star$ was proved to be associative in Subsection~\ref{ssec-gengaugegroupoid}.
Moreover, for any HS morphism $P$, the inverse in $\calH$ of the division map $\phi_P$ is a HS generalized gauge transformation on $P$ by Proposition~\ref{prop-proddivhils}, which is obviously the image of the HS gauge transformation $\iota_{\calH}\circ \varepsilon_P$ of $P$.
Any HS generalized gauge transformation $K$ between any two object of $\mathsf{HS}_{\calG,\calH}$ was proved to be invertible; clearly, this translates into the fact that the two possible products $\star$ between $K$ and its inverse equal the inverse in $\calH$ of the division maps of the HS morphisms to which $K$ is attached.

Therefore, given two groupoids $\calG$ and $\calH$, it makes sense to define the {\em groupoid of HS generalized gauge transformations $C^{\infty,\calH^2}_{\calG}$} by the following data:
\begin{itemize}
\item[i)] the objects of $C^{\infty,\calH^2}_{\calG}$ are the objects of the category $\mathsf{HS}_{\calG,\calH}$. 
\item[ii)] Given any two objects $P_1$, $P_2$ of $C^{\infty,\calH^2}_{\calG}$, the set of arrows is the set $C^{\infty}_{\calG}\!\left(P_1\odot P_2,\calH\right)^{\calH^2}$ of HS generalized gauge transformations between $P_1$ and $P_2$.  
\item[iii)] The target $t$, source $s$ and unit map $\iota$ are defined as for the groupoid of generalized gauge transformations introduced in subsection~\ref{ssec-gengaugegroupoid}; the product for composable HS generalized gauge transformations is set to be $\star$, as in Subsection~\ref{ssec-gengaugegroupoid}.
\end{itemize}
With these data, and using the same arguments of Subsection~\ref{ssec-gengaugegroupoid}, it follows
\begin{Thm}
The $6$-tuple $\left(C^{\infty,\calH^2}_{\calG},\mathsf{HS}_{\calG,\calH},s,t,i,j\right)$, for any two groupoids $\calG$ and $\calH$, is an abstract groupoid in the sense of Definition~\ref{def-groupoid}; it is obviously isomorphic to the abstract groupoid $\mathsf{HS}_{\calG,\calH}$, a category whose morphisms are all invertible, by Theorem~\ref{thm-gengaugehils}.
For any object $P$ of $\mathsf{HS}_{\calG,\calH}$, the isotropy group $C^{\infty,\calH^2}_{\calG,P}$ is isomorphic to the HS gauge group of $P$, $C^{\infty}_{\calG}\!\left(P,\calH\right)^{\calH}$.
\end{Thm}

\end{document}